\documentclass[10pt,a4paper]{elsarticle}
\usepackage[utf8]{inputenc}
\usepackage{amsmath}
\usepackage{amsfonts}
\usepackage{amssymb}
\usepackage{graphicx}
\usepackage{subcaption}
\usepackage{color, colortbl}
\usepackage[first=0,last=9]{lcg}

\definecolor{Gray}{gray}{0.9}

\begin{document}

\begin{frontmatter}

\title{A universal centred high-order method based on implicit Taylor series expansion with fast second order evolution of spatial derivatives
}

\author{Gino I. Montecinos\corref{cor1}}
\ead{gino.montecinos@uaysen.cl}

\cortext[cor1]{Corresponding author}
\address[myuni]{Universidad de Ays\'en}



\begin{abstract}
In this paper, a centred universal high-order finite volume method for solving hyperbolic balance laws is presented. 
The scheme belongs to the family of ADER methods where the Generalized Riemann Problems (GRP) is a building block. The solution to these problems is carried through an implicit Taylor series expansion, which allows the scheme to works very well for stiff source terms. A von Neumann stability analysis is carried out to investigate the range of CFL values for which stability and accuracy are balanced. The scheme implements a centred, low dissipation approach for dealing with the advective part of the system which profits from small CFL values. Numerical tests demonstrate that the present scheme can solve, efficiently, hyperbolic balance laws in both conservative and non-conservative form as well. An empirical convergence rate assessment shows that the expected theoretical orders of accuracy are achieved up to the fifth order. 

\end{abstract}


\end{frontmatter}

\section{Introduction}

Hyperbolic balance laws play a crucial role in describing phenomena in several fields and applications. These systems are characterized by the presence of source terms, where several issues arise under both theoretical and numerical points of view. From the theoretical point of the view, the existence of solutions is one of the most challenging issues, \cite{Kruzkov:1970a, Liu1987, Liu:1981b}. Already in the scalar case, for some particular source terms the solution can blow up in finite time, also the gradient of the solution can suffer from this phenomena, \cite{Bressan:2000a}. However, even if the analytical solutions do not depict this behavior, the numerical approximations may suffer from the presence of source terms, a good example is a class of problems called stiff, where the reactive characteristic speed associated with the source terms are faster than the characteristic speeds associated to the convective part of hyperbolic problems.   Stiff source terms can appear artificially in numerical solutions of hyperbolic laws even if the source term is non-stiff or it is not present, as occurs for instance in \cite{Jin:1995a, Kawashima:1987a, Bouchut:1999a} or in the case of hyperbolic reformulations of parabolic problems \cite{Montecinos:2014a, Montecinos:2014b}. So, independently on the nature of the source term, the construction of numerical schemes able to efficiently obtain solutions of hyperbolic balance laws is an issue of current interest.

On the other hand, the hyperbolic problems always can be expressed in the so called quasilinear form, \cite{Bressan:2000a,Toro:2009a}. If there exists some flux function such that its Jacobian, corresponds to the matrix involved on the quasilinear form, the system is said to be conservative. On the other hand, if a such function there is not exist, then the system is said to be non-conservative, \cite{DalMaso:1995a, CastroMJ:2006a}. Non-conservative schemes are another important issues under the numerical point of view, \cite{DalMaso:1995a,Abgrall:2010a, CastroMJ:2006a, Pares:2006a}. Non conservatives systems are characterized by containing non-conservative products in the sense of \cite{DalMaso:1995a} and the concept of numerical fluxes makes no sense and has to be extended to the non-conservative counterpart, the numerical fluctuations or increments. Since conservation laws can be expressed in quasilinear form through the Jacobian matrix, a desirable property of numerical schemes for solving non-conservative systems is that these are able to act as a conservative method on those cases. It is well known that some path-conservative types schemes \cite{DalMaso:1995a,Pares:2006a}, have this duality. A numerical scheme is called universal, in the sense of \cite{Balsara:2018b} if this can act as a conservative one, without any modification, if there exists a flux function whose Jacobian is used to form quasilinear systems.

The ADER ({\bf A}rbitrary Accuracy {\bf DER}ivative Riemann
problem) schemes , \cite{Millington:2001a, Toro:2002a, Titarev:2002a, Toro:2006b} is a family of high-order finite volume scheme which works well for general hyperbolic balance laws. These can incorporate both the technology for dealing with stiff source terms \cite{Dumbser:2008a, Balsara:2009a, Dumbser:2013a, Boscheri:2014c, Montecinos:2014a, Toro:2015a} and the universal property to deal with non-conservative systems using the path-conservative approach \cite{Dumbser:2009b}. These schemes require the solution of the so called Generalized Riemann Problems (GRP), whose solutions can be found analytically or in an approximate form. These are then used to evaluate the source terms and the numerical flux, in the case of conservative problems, or the so called numerical fluctuations, in the case of non-conservative problems.  

Generalized Riemann Problems (GRP) are initial value problems and these are extensions of the classical one by the incorporation of two elements; i) the presence of the source term and; a piecewise smooth initial condition. Classical Riemann problems are initial value problems in which the PDE is a homogeneous equation and the initial condition is a piecewise constant function, see \cite{Toro:2009a, Toro:2016a} for further details on Riemann problems as a building block for finite volume methods.

To solve GRP there are two approaches, the first one due to Toro and Titarev named here the TT solver \cite{Toro:2002a, Toro:2006b}, where the solution is expressed in terms of a Taylor series expansion in time. The time derivatives are completely expressed in terms of state equations and spatial derivatives of them. This is achieved by a systematic use of the governing equation, this strategy corresponds to the Cauchy-Kowalewskaya procedure and consists of constructing a set of functional representing the time derivatives, whose arguments are the solution to a classical Riemann problem defined for the leading term and a sequence of linearized Riemann problems for the spatial derivatives.

The second approach for solving GRP is due to an adaptation of the second order method of Harten et al \cite{Harten:1987b}, presented by Castro and Toro in \cite{Castro:2008a}, named here the HEOC solver.  In this solver, local predictors of the PDE are obtained within the computation cell. Then local classical Riemann problems are constructed on the cell interfaces from the homogeneous part of the governing equation (zeroing the source term) and the piecewise constant initial condition is obtained from the predictors evaluated at precise locations on the interface. From the solution of these classical Riemann problems, it is evaluated the numerical flux or the numerical fluctuation, for conservative and non-conservative problems respectively.  The numerical source is straightforwardly evaluated from the predictor within the computational cell. A common factor in these two strategies is the reconstruction procedure. See \cite{Toro:2008a, Montecinos:2012a} for further details about GRP solvers.

The predictor can be computed by using the Taylor series expansion as in the TT solver \cite{Toro:2002a} and discontinuous Galerkin approach \cite{Dumbser:2008a}, as well. The approach based on Galerkin formulations works well for stiff source terms but it requires the conventional matrix inversion of the Finite element framework. The strategy using Taylor series expansions may be more efficient but in principle does not work for stiff source terms. Furthermore, the Cauchy-Kowalewskaya procedure becomes cumbersome for high order cases or complex systems.  However, Toro and Montecinos \cite{Toro:2015a}  have shown how implicit Taylor series can be adapted for dealing with stiff source terms. Recently, in \cite{Montecinos:2020a} a simplified Cauchy-Kowalewsky which only requires the Jacobian matrices and a strategy for evaluation their derivatives are enough to generate all the functional in a very efficient way.

In \cite{Toro:2015a} implicit Taylor series expansions are used to express the solution of the GRP, this strategy requires the implicit evolution of the spatial derivatives as well. The Taylor series for both the state and all spatial derivatives are simultaneously incorporated into a very large algebraic system, which is solved by using fixed point iteration procedures.  The size of the system involves the number of unknowns and the order of accuracy as well. The complexity of such an approach increases exponentially with respect to the order of accuracy. So the efficiency of this approach is penalized as the order increases.

Despite a new strategy for approximating the Cauchy-Kowalewskaya functional is available in \cite{Montecinos:2020a}, the existence of the closed forms given by the conventional Cauchy-Kowalewskaya functional should be as efficient as the simplified version. For this reason, in this work a practical strategy for implementing the conventional procedure for general hyperbolic balance laws, based on the syntaxes of the well known, GNU computer algebra system {\it Maxima}, \cite{Ochsner:2019a}, is reported and so this is the approach implemented in this paper.

In this work a variation of the scheme in \cite{Toro:2015a} is presented, this uses implicit Taylor series for the data and the conventional Cauchy-Kowalewskaya procedure. However, for evolving spatial derivatives we use a second order evolution by using a suitable linearization of the governing equation. This simplification generates a non-linear algebraic system, but this depicts a very nice structure that allows us to design an inversion matrix procedure, where a matrix of size proportional to the number of unknowns, for any order of accuracy, needs to be evolved. This yields an improvement of the efficiency of the scheme. Even if, from a von Neumann stability analysis a small CFL coefficient is required, the global efficiency is superior to that of the approach in \cite{Toro:2015a} as the accuracy increases.  This is both analytically and theoretical demonstrated.

Recently, a low-dissipation, first order, centred scheme presented by Toro et al \cite{Toro:2020a} profits from low values of CFL.  It is well known that centred schemes fail in the case of linearly degenerated, intermediate, characteristic field. However, these are suitable methods if the complete solution of the Riemann problems is not available, see \cite{Toro:2016a} for further discussion about the role of the solution of Riemann problems on the computation of numerical fluxes.  As shown in \cite{Toro:2016a}  and \cite{Toro:2020a}, the FORCE scheme, the centred method on which the scheme implemented here is based on, is less diffusive than the HLL solver and the Rusanov one, respectively.  In \cite{Toro:2020a} is illustrated the advantages of this novel scheme on small CFL regimes in terms of the numerical dissipation. It is demonstrated, for the linear advection equation that dissipation of the scheme can be comparable with that of the Godunov method, the lowest one reported in the literature. Despite the method is devoted to first order scheme, this is a suitable approach to be explored on the present method.

The aim of this article is to present a numerical scheme for balance laws, centred, universal and based on the implicit Taylor series expansion and the conventional Cauchy-Kowalewskaya procedure which makes this to be able to deal with stiff source terms achieving arbitrary high-order of accuracy.

This document is organized as follows. In section \ref{sec:the-method} the framework of the numerical approach is presented. In section \ref{sec:the-predictor} the predictor for the GRP solver is implemented. In section \ref{sec:theoretical-analysis}, the von Neumann stability analysis and an estimation of the efficiency are carried out. In section \ref{sec:num-results}, numerical results are presented. In section \ref{sec:conclusions} the conclusions are drawn.

\section{The method}\label{sec:the-method}
In this document, it is proposed an alternative formulation to the GRP-solver based on implicit Taylor series expansion and Cauchy-Kowalewskaya procedure proposed by Toro and Montecinos (2016). The GRP of interest in this work is to be used as a building block in the construction of high order scheme for solving general hyperbolic balance laws of the form
\begin{eqnarray}
\label{eq:gov-equation}
\begin{array}{c}
\partial_t \mathbf{Q} (x,t) + \partial_t \mathbf{F}(\mathbf{Q}(x,t))  +  \mathbf{B}(\mathbf{Q}(x,t)) \partial_x \mathbf{Q}(x,t) = \mathbf{S}(\mathbf{Q}(x,t))\;.
\end{array}
\end{eqnarray}
We can also express these problems as 
\begin{eqnarray}
\label{eq:gov-equation-non-cons}
\begin{array}{c}
\partial_t \mathbf{Q} (x,t) + \mathbf{A}(\mathbf{Q}(x,t)) \partial_x \mathbf{Q}(x,t) = \mathbf{S}(\mathbf{Q}(x,t))\;,
\end{array}
\end{eqnarray}
with $\mathbf{A}(\mathbf{Q}) = \frac{\partial \mathbf{F}(\mathbf{Q})}{\partial \mathbf{Q}} + \mathbf{B}(\mathbf{Q})\;.$
Therefore, we can apply the scheme presented by Toro et al. \cite{Toro:2020a}, which is based on the general form of a path-conservative scheme
\begin{eqnarray}
\label{eq:FV-non-cons}
\begin{array}{c}

\mathbf{Q}_i^{n+1}  = \mathbf{Q}_i^n - \frac{\Delta t}{ \Delta x} ( \mathbf{D}_{i-\frac{1}{2}}^{+} +  \mathbf{D}_{i+\frac{1}{2}}^{-} )   + \Delta t ( \mathbf{S}_{i} - \mathbf{A}_i^n)\;,
\end{array}
\end{eqnarray}
where $\Delta t = t^{n+1} - t^n$, $\Delta x = x_{i+\frac{1}{2}} - x_{i-\frac{1}{2}} $ and
\begin{eqnarray}
\begin{array}{c}
\displaystyle
\mathbf{D}_{i + \frac{1}{2}}^{\pm} 
= \frac{ 1}{ \Delta t} \int_{t^n }^{t^{n+1}}
\mathbf{A}_{i + \frac{1}{2} }^{\pm}(  \mathbf{Q}_{i}( x_{i+\frac{1}{2}},\tau), \mathbf{Q}_{i+1}( x_{i+\frac{1}{2}},\tau)) 
\cdot ( \mathbf{Q}_{i+1}( x_{i+\frac{1}{2}},\tau)-  \mathbf{Q}_{i}( x_{i+\frac{1}{2}},\tau))
d \tau
\end{array}
\end{eqnarray}
and 
\begin{eqnarray}
\begin{array}{c}
\displaystyle
\mathbf{S}_{i} 
= \frac{ 1}{ \Delta x \Delta t} 
\int_{x_{i-\frac{1}{2}}}^{x_{i+\frac{1}{2}}}\int_{t^n }^{t^{n+1}}
\mathbf{S}( \mathbf{Q}_{i}(  \mathbf{Q}_{i}(\xi,\tau) )  d \tau d\xi

\;,
\end{array}
\end{eqnarray}
\begin{eqnarray}
\label{eq:integration-Ax}
\begin{array}{c}
\displaystyle
\mathbf{A}_{i}^n 
= \frac{ 1}{ \Delta x \Delta t} 
\int_{x_{i-\frac{1}{2}}}^{x_{i+\frac{1}{2}}}\int_{t^n }^{t^{n+1}}
\mathbf{A}( \mathbf{Q}_{i}( \xi,\tau) ) \cdot \partial_x \mathbf{Q}_{i} ( \xi,\tau)  d \tau d\xi

\;,
\end{array}
\end{eqnarray}
Here, $\mathbf{Q}_i(x,t)$ is the predictor within the space-time cell $[x_{i-\frac{1}{2}}, x_{i+\frac{1}{2}}] \times [t^n, t^{n+1} ]$. The expression $\mathbf{A}_i^n$   appears only in high order formulations, this is zero in the first order case. Here, $\mathbf{A}_{i+\frac{1}{2}}^{\pm}$  can be seen as the two state matrix function
\begin{eqnarray}
\begin{array}{c}
\displaystyle
\mathbf{A}_{i+\frac{1}{2}}^{\pm}(\mathbf{Q}_L, \mathbf{Q}_R) = \frac{1}{2} \tilde{\mathbf{A}}( \mathbf{Q}_L, \mathbf{Q}_R)
\pm \frac{\alpha \Delta t}{4 \Delta x} \biggl[ \tilde{\mathbf{A}}( \mathbf{Q}_L, \mathbf{Q}_R)^2  + \biggl( \frac{ \Delta x}{\alpha \Delta t} \biggr)^2 \mathbf{I}\biggr] \;,
\end{array}
\end{eqnarray}
where $\mathbf{I}$ is the identity matrix and 
\begin{eqnarray}
\begin{array}{c}
\displaystyle
\tilde{\mathbf{A}}( \mathbf{Q}_L, \mathbf{Q}_R)

= \int_0^1 \mathbf{A}( \Psi(s; \mathbf{Q}_L, \mathbf{Q}_R) ) ds  \;.
 
\end{array}
\end{eqnarray}
This expression is obtained from the use of the segment path $\Psi (s, \mathbf{Q}_L, \mathbf{Q}_R) = \mathbf{Q}_L + s \cdot ( \mathbf{Q}_R - \mathbf{Q}_L)$. From a suitable quadrature rule consisting of $n_{GP}$ points $\xi_j$ and corresponding wights $\omega_j$, this integral can be computed as 
\begin{eqnarray}
\begin{array}{c}
\displaystyle
\tilde{\mathbf{A}}( \mathbf{Q}_L, \mathbf{Q}_R)
 = 
 \sum_{j=1}^{nGP} \omega_j \mathbf{A}( \Psi( \xi_j; \mathbf{Q}_L, \mathbf{Q}_R) ) \;.
 
\end{array}
\end{eqnarray}
Notice that this is the high-order counterpart of the scheme in \cite{Canestrelli:2009a} but adapted to the low-dissipation order scheme \cite{Toro:2020a}.

For the first order case, has been proved in \cite{Toro:2020a} that the scheme depicts a numerical diffusion which is comparable with that of Godunov, one of the lowest dissipative schemes available in the literature.

\section{The predictor step}\label{sec:the-predictor}
In this section we are going to consider the procedure to get the predictor $\mathbf{Q}_{i}(x,t)$ into the computation space-time cell $ I^n_i = [x_{i-\frac{1}{2}}, x_{i+\frac{1}{2}}]\times [t^{n} ,t^{n+1}]$ in terms of the implicit Taylor series expansion, as in Toro and Montecinos \cite{Toro:2015a}. However, for the sake of simplicity we are going to omit the subscript $i$. Then, we are going to apply the change of variables $ x = x_{i- \frac{1}{2}} + \xi \Delta x$ and $t = t^n + \tau \Delta t$ in order to translate $I^n_i$ into $[0,1]^2$. So the predictor at any location $(\xi, \tau)$ in $[0,1]^2$ is obtained as 
\begin{eqnarray}
\label{eq:taylor-1}
\begin{array}{c}
\displaystyle
\mathbf{Q}(\xi, \tau) = \mathbf{Q}(\xi, 0 ) - \sum_{ k = 1}^{ M }  \frac{ (- \tau)^{ k} }{ k!} \partial_t^{(k)} \mathbf{Q}( \xi, \tau) \;.  
\end{array}
\end{eqnarray} 
For evaluating the time derivative, we use the  Cauchy-Kowalewskaya procedure. We can use the recently simplified version of the Cauchy-Kowalewskaya procedure in \cite{Montecinos:2020a}. However, once the Cauchy-Kowalewskaya functionals are available the resulting scheme becomes a very efficient one. Since, an external tool for generating these functionals is reported in the appendix \ref{sec:ck-generator}, where the cumbersome procedure is dramatically simplified. We insist on using the conventional Cauchy-Kowalewskaya procedure.   In any option one choices, the Taylor series expansion (\ref{eq:taylor-1}) can be written as     
\begin{eqnarray}
\label{eq:taylor-and-CK}
\begin{array}{c}
\displaystyle
\mathbf{Q}(\xi, \tau) = \mathbf{Q}(\xi, 0 ) - \sum_{ k = 1}^{ M }  \frac{ (- \tau)^{ k} }{ k!} 
\mathbf{G}^{(k)}( \mathbf{Q}( \xi, \tau), \partial_x \mathbf{Q}( \xi, \tau), ..., \partial_x^{(k)} \mathbf{Q}( \xi, \tau) ) ,   \;.  
\end{array}
\end{eqnarray}
where $ \mathbf{G}^{(k)}$ corresponds to the Cauchy-Kowalewskaya functional for the $k$th time derivative. Since, in \cite{Toro:2015a, Montecinos:2012a} has been shown that for stability purposes, the spatial derivatives must be considered at the same time  $\tau$ than the state. So, (\ref{eq:taylor-and-CK}) requires the evolution in time of the spatial derivatives of all orders, up to $ M$. 

In this paper, we are going to propose a strategy that only will require the inversion of two $m \times m$ matrices.  The strategy is based on a second order variation in time of the spatial derivatives. To stem the procedure we need to make some simplification. 
To simplify the notation, let us going to denote by $\mathbf{Q}^*$ the solution state at $(\xi, \tau)$, that is, the constant state obtained once the solution has already been evaluated at those coordinates. The first simplification assumes that near the solution state the partial differential equation (\ref{eq:gov-equation}) behaves as 
\begin{eqnarray}
\label{eq:gov-equation-linearised}
\begin{array}{c}
\partial_t \mathbf{Q} (x,t) + \mathbf{A}^* \partial_t \mathbf{Q}(x,t) = \mathbf{S}(\mathbf{Q}(x,t))\;,
\end{array}
\end{eqnarray}
where $ \mathbf{A}^* = \mathbf{A}(\mathbf{Q}^*) = \frac{\partial \mathbf{F}(\mathbf{Q}^*)}{\partial \mathbf{Q}}$, that is, we can linearise the evolution around the solution state.
We differentiate (\ref{eq:gov-equation-linearised}) with respect to the $x$ variable and obtain   
\begin{eqnarray}
\label{eq:gov-equation-dx}
\begin{array}{c}
\partial_t (\partial_x \mathbf{Q}(x,t)) + \mathbf{A}^* \partial_t ( \partial_x \mathbf{Q}(x,t) ) 
= \mathbf{B}(\mathbf{Q}(x,t)) \partial_x \mathbf{Q}(x,t)\;,
\end{array}
\end{eqnarray}
where $\mathbf{B}$ is the Jacobian matrix of $\mathbf{S}$ with respect to $\mathbf{Q}$. Here, again, as the system  (\ref{eq:gov-equation}) near the state $ \mathbf{Q}^*$ behaves as (\ref{eq:gov-equation-linearised}), we can assume that (\ref{eq:gov-equation-dx}) behaves as   
\begin{eqnarray}
\label{eq:gov-equation-dx-linearised}
\begin{array}{c}
\partial_t (\partial_x \mathbf{Q}(x,t)) + \mathbf{A}^* \partial_t ( \partial_x \mathbf{Q}(x,t) ) 
= \mathbf{B}^* \partial_x \mathbf{Q}(x,t)\;,
\end{array}
\end{eqnarray}
with $\mathbf{B}^* = \mathbf{B}(\mathbf{Q}^*) $.
The relevant assumption of this linearized equation is that we can differential this system ( assuming that the system is regular enough) to obtain
\begin{eqnarray}
\label{eq:gov-equation-dx-k-linearised}
\begin{array}{c}
\partial_t (\partial_x^{(k)} \mathbf{Q}(x,t)) + \mathbf{A}^* \partial_t ( \partial_x^{(k)} \mathbf{Q}(x,t) ) 
= \mathbf{B}^* \partial_x^{(k)} \mathbf{Q}(x,t)\;.
\end{array}
\end{eqnarray}
From this linearised system, we can assert that a second order evolution of spatial derivatives has the form
\begin{eqnarray}
\label{eq:taylor-dx-k}
\begin{array}{c}
\displaystyle
\partial_x^{(k)}\mathbf{Q}(\xi, \tau) = \partial_x^{(k)} \mathbf{Q}(\xi, 0 ) + \tau \left(

\mathbf{B}^* \partial_x^{(k)} \mathbf{Q}(\xi, \tau) - \mathbf{A}^* \partial_x^{(k+1)}\mathbf{Q}(\xi, \tau)

 \right)\;.  
\end{array}
\end{eqnarray}
This is true for any $k$, however, for practical implementations, we have only access to $M$ spatial derivatives. So, we are going to neglect spatial derivatives of orders higher than $M$. We do not have access to $M$ spatial derivatives because if we assume initial conditions to be polynomial of orders at most $M$ then, derivatives higher than this order are zero.   For the evolution of the $M$th spatial derivative, we propose
\begin{eqnarray}
\label{eq:taylor-dx-M}
\begin{array}{c}
\displaystyle
\partial_x^{(M )}\mathbf{Q}(\xi, \tau) = \partial_x^{(M )} \mathbf{Q}(\xi, 0 ) + \tau \left(

\mathbf{B}^* \partial_x^{(M )} \mathbf{Q }(\xi, \tau) 

 \right)\;.  
\end{array}
\end{eqnarray}

On the other hand, to obtain $\mathbf{Q}^*$ from (\ref{eq:taylor-and-CK}) we require $\partial_t^{(k)}\mathbf{Q}$, for $k=1,...,M$.  So, equations (\ref{eq:taylor-and-CK}), (\ref{eq:taylor-dx-k}) and (\ref{eq:taylor-dx-M})  can be set into a non-linear algebraic equation for the state and their spatial derivatives. We can solve the full resulting algebraic equation. However, this can be as expensive as resolving the original system in Toro and Montecinos \cite{Toro:2015a}. Instead of the full system we proposed a procedure which only involves the inversion of two  $m\times m$ small matrices. To introduce the idea, let us define vectors $\mathbf{D}_k $ as $ \mathbf{D}_k = \partial_x^{(k)} \mathbf{Q}$, the convention $ \partial_x^{(0)} \mathbf{Q} = \mathbf{Q}$ has  been used.
Once this definition has been done, the  equations (\ref{eq:taylor-and-CK}), (\ref{eq:taylor-dx-k}) and (\ref{eq:taylor-dx-M})  can be put together as
\begin{eqnarray}
\label{eq:system-reormulated}
\begin{array}{lcl}

\displaystyle

\mathbf{D}_0 &=& \mathbf{Q}(\xi, 0 ) - \sum_{ k = 1}^{ M }  \frac{ (- \tau)^{ k} }{ k!} 
\mathbf{G}^{(k)}( \mathbf{D}_0,  \mathbf{D}_1, ...,  \mathbf{D}_k )    \;,  \\

&\vdots &  \\
\displaystyle

\mathbf{D}_k &=& \partial_x^{(k)} \mathbf{Q}(\xi, 0 ) + \tau \left(

\mathbf{B}(\mathbf{D}_0) \mathbf{D}_k - \mathbf{A}(\mathbf{D}_0) \mathbf{D}_{k+1}

 \right)  \; \\
 
 &\vdots &  \\
\displaystyle
\mathbf{D}_M &=& \partial_x^{(M )} \mathbf{Q}(\xi, 0 ) + \tau \left(

\mathbf{B}(\mathbf{D}_0) \mathbf{D}_M  \right) \;.

\end{array}
\end{eqnarray}
We propose the following nested fixed point procedure
\begin{eqnarray}
\label{eq:nested-procedure}
\begin{array}{lcl}

\displaystyle

\mathbf{D}_0^{r+1} &=& \mathbf{Q}(\xi, 0 ) - \sum_{ k = 1}^{ M }  \frac{ (- \tau)^{ k} }{ k!} 
\mathbf{G}^{(k)}( \mathbf{D}_0^{r+1},  \mathbf{D}_1^{r}, ...,  \mathbf{D}_k^{r} )    \;,  \\

&\vdots &  \\
\displaystyle

\mathbf{D}_k^{r+1} &=& \partial_x^{(k)} \mathbf{Q}(\xi, 0 ) + \tau \left(

\mathbf{B}(\mathbf{D}_0^{r}) \mathbf{D}_k^{r+1} - \mathbf{A}(\mathbf{D}_0^{r}) \mathbf{D}_{k+1}^{r+1}\;,

 \right)  \; \\
 
 &\vdots &  \\
\displaystyle
\mathbf{D}_M^{r+1} &=& \partial_x^{(M )} \mathbf{Q}(\xi, 0 ) + \tau \left(

\mathbf{B}(\mathbf{D}_0^{r}) \mathbf{D}_M^{r+1}  \right) \;.

\end{array}
\end{eqnarray}
where $r$ is an iteration index.  
As in the case of Toro and Montecinos \cite{Toro:2015a}, $\partial_x^{(k)} \mathbf{Q}(x,0) = \mathbf{P}^{(k)}(x)$, where
 $\mathbf{P}(x)$ is the reconstruction polynomial defined into the computational cell where the predictor is being obtained.
Once we implement a descent step type scheme. The Jacobian matrix of the algebraic systems for $\mathbf{D}_k$, with $ k=1,...,M$ depicts a very nice structure which allows us to build a very efficient inversion matrix procedure, where only one $m\times m$ matrix is required. The evolution of $\mathbf{D}_0$ requires the determination of the Cauchy-Kowalewskaya procedure as well as the derivatives of these functionals with respect to $\mathbf{D}_0$. In the appendix \ref{sec:ck-generator} a practical script to generate the Cauchy-Kowaleskaya and the corresponding Jacobian functional of the algebraic system resulting from the implicit Taylor series, is reported. The methodology is general enough to include a general hyperbolic system.

Up to here, the description for obtaining the predictor $\mathbf{Q}(\xi, \tau)$ has been presented. And thus, the expressions for (\ref{eq:FV-non-cons}) can be obtained. Notice that (\ref{eq:integration-Ax}) requires not only the predictor but also the derivative of this. Since the value $  \mathbf{ Q}(\xi, \tau)$ at each quadrature point in space and time is known, we use an interpolation polynomial of $\mathbf{Q} $ as in \cite{Montecinos:2020a}, that is, at each quadrature point $\tau_j $ in time we use the set $ \{ \mathbf{Q}(\xi_l, \tau_j)\}_{l}$ to carry out a Lagrange interpolation.  Then, we use the derivative of this polynomial to approximate $\partial_x \mathbf{Q}$. On the other hand, from the solution of (\ref{eq:system-reormulated}), $\mathbf{D}_k$ is a good candidate, however, this is only a second order approximation, which is enough for the implicit Taylor series (\ref{eq:taylor-1}) but it is not so for evaluating (\ref{eq:integration-Ax}).

\section{Theoretical analysis of the present approach}\label{sec:theoretical-analysis}
In this section, the von Neumann stability analysis and the analytical estimation of the number of operations and how it affects the global efficiency of the present scheme, are investigated.

\subsection{The von-Neumann stability analysis}\label{sec:von-neumann}
In this section, it is presented a new strategy for computing the von Neumann stability analysis for high-order ADER schemes using the WENO reconstruction procedure and applied to the lineal advection reaction equation
\begin{eqnarray}
\begin{array}{c}
\partial_t q + \lambda \partial_x q = \beta q \;,
\end{array}
\end{eqnarray}
with $\lambda > 0$ and $\beta \leq 0$, constant values. The present approach will be applied to the ADER method using the HEOC approach. However, the similar study can be applied also to the original ADER scheme of Toro and Titarev, \cite{Toro:2001c, Toro:2002a, Titarev:2002a}.  Notice that for the case of the linear advection conservative and non-conservative formulations as presented in this work, are equivalent. Therefore, we are going to be interested on conservative methods 
\begin{eqnarray}
\begin{array}{c}
q_i^{n+1} 
=
q_i^{n} - \frac{ \Delta t}{ \Delta x} \biggl(  
f_{i+\frac{1}{2}} - f_{i-\frac{1}{2}}
\biggr) 
+\Delta t \cdot s_i \;,
\end{array}
\end{eqnarray}
with general numerical flux functions. That is, the numerical flux is obtained as 
\begin{eqnarray}
\begin{array}{c}
f_{i+\frac{1}{2}}
= \sum_{ u = 1 }^{ n_{GP} } \omega_u \cdot F_{RP} ( q_{i}(1,\eta_u), q_{i+1}(0,\eta_u)) \;,
\end{array}
\end{eqnarray}
where $(\eta_u, \omega_m)$ represents the pair of quadrature points and quadrature weights, respectively. Here, we use the Gauss-Lobatto quadrature rule. The expression $ F_{RP }( q_L, q_R ) \;,$ represents any numerical flux function of two states $q_L$ and $q_R$, resulting from an approximate Riemann solver of classical Riemann problems.  
The source term is computed as
\begin{eqnarray}
\begin{array}{c}

\displaystyle

s_{ i}
=
\beta \cdot \sum_{ v = 1 }^{ n_{GP} } \sum_{ u = 1 }^{ n_{GP} } \omega_u \cdot \omega_v \cdot q_i(\eta_v, \eta_u )

 \;.
\end{array}
\end{eqnarray}

The expression $q_{i}(\xi, \tau)$ represents the predictor within the computational cell $[x_{i-\frac{1}{2}},x_{i+\frac{1}{2}}]\times [t^n, t^{n+1}]$ in terms of local coordinates $(\xi,\tau)$ on $[0,1]$.  Let us introduce $ c = \lambda \cdot \frac{ \Delta t}{ \Delta x}$ and  $ r = \beta \cdot \Delta t$, then we introduce $q_i^n = A^n \cdot e^{(I \cdot \theta \cdot i)}$ with $ \theta $ the phase angle with $ \theta \in [0, 2 \pi)$, $I^2 = -1$ the unity in the complex numbers and $A$ is a complex amplitude. So, finite volume formula provides the following function for the amplitude 
\begin{eqnarray}
\label{eq:amplitude-ADER}
\begin{array}{c}
A(\theta, c, r) := 1  - A^{-n} \cdot e^{-I \cdot \theta \cdot i} \cdot\biggl( c \cdot ( \hat{f}_{i+\frac{1}{2}}(\theta, c, r) - \hat{f}_{i-\frac{1}{2}}(\theta, c, r)  ) + r \cdot  \hat{s}_i(\theta, c, r)  \biggr)\;,
\end{array}
\end{eqnarray}
where
\begin{eqnarray}
\begin{array}{c}

\displaystyle

\hat{f}_{i+\frac{1}{2}}(\theta, c, r)
=
\sum_{ m = 1 }^{ n_{GP} } \omega_m \cdot 

F_{RP}( q_i(\theta, c, r, 1, \eta_m ), q_{i+1}(\theta, c, r, 0, \eta_m ) )

\end{array}
\end{eqnarray}
and

\begin{eqnarray}
\begin{array}{c}

\displaystyle

\hat{s}_{ i} (\theta, c, r)
=
\sum_{ v = 1 }^{ n_{GP} } \sum_{ m = 1 }^{ n_{GP} } \omega_m \cdot \omega_v \cdot q_i(\theta, c, r, \eta_v, \eta_m )

 \;.
\end{array}
\end{eqnarray}

A scheme is stable for $c$ and $r$ if $|A| \leq 1$ for any $\theta \in [0, 2 \pi)$, here $|A| = \sqrt{ Im(A)^2 + Re(A)}$ is the module of the complex number, $Im(A)$ and $Re(A)$  are the imaginary and real parts, respectively.

\subsubsection{The explicit predictor}
In this section, we are going to analyse the high order ADER-HEOC. As illustrated in \cite{Castro:2008a,Montecinos:2012a} and detailed in the section \ref{sec:the-method} of the current paper, the numerical fluxes involve the solution of one classical Riemann problem which is constructed from local predictors within computational cells. This predictor is expressed in terms of the Taylor series expansion
\begin{eqnarray}
\begin{array}{c}
q_{i}(\xi,\tau)
=

w_i(\xi) + \sum_{k=1}^M  \frac{ ( \tau)^k }{ k!} \partial_t^{(k)} q_i( \xi, 0)
\;,
\end{array}
\end{eqnarray}
where $ \partial_t^{(k)} q_i( \xi, 0) = G^{(k)} (q_i(\xi,0), \partial_x q_i(\xi,0),...,\partial_x^{(k)} q_i(\xi,0))$ is the  Cauchy-Kowalewskaya functional, which for the case of the linear advection reaction equation has the explicit form
\begin{eqnarray}
\begin{array}{c}
\displaystyle
 G^{(k)} (q_i, \partial_x q_i,...,\partial_x^{(k)} q_i)  := \sum_{l=0}^k 

\sum_{j = 1}^l \left(
\begin{array}{c}

l \\
j

\end{array}
\right) 
\partial_x^{(j)} q_i(\xi,\tau) \cdot (-\lambda)^{j} \cdot (\beta)^{l-j}\;.
\end{array}
\end{eqnarray}
Instead of $q_i(\xi,0)$, we use the reconstruction polynomial $w_i(\xi)$, so the Cauchy-Kowalewskaya functional takes the form
\begin{eqnarray}
\begin{array}{c}
\displaystyle
 G^{(k)} = 
 
\sum_{l=0}^k 
\sum_{j = 1}^l \left(
\begin{array}{c}

l \\
j

\end{array}
\right) 
\frac{ d^{j} }{dx^{j}} w_i(\xi) \cdot (-\lambda)^{j} \cdot (\beta)^{l-j}

\;.
\end{array}
\end{eqnarray}
Therefore, the predictor can be written as
\begin{eqnarray}
\begin{array}{lcl}

\displaystyle

q_{i}(\xi,\tau)
=
  w_i(\xi)

+ \sum_{k=1}^M  \frac{ ( \tau)^k }{ k!} 

\sum_{l=0}^k 
\sum_{j = 1}^l \left(
\begin{array}{c}

l \\
j

\end{array}
\right) 
\frac{ d^{j} }{dx^{j}} w_i( \xi) \cdot (-\lambda)^{j} \cdot (\beta)^{l-j}\;.

\end{array}
\end{eqnarray}

Now, we carry out the von Neumann analysis introducing $ q^n_i = A^n e^{  I \cdot  \theta \cdot i }$, where $ \theta$ is the phase angle, $ I^2 = -1 $ is the imaginary unity and $ A $ is the amplitude, a complex number.

In \ref{sec:reconstruction}, a review of the reconstruction procedure is carried out. It is also shown that $w_i(\xi)$ can be expressed as
\begin{eqnarray}
\begin{array}{c}
\displaystyle
w_i(\theta,\xi) =  \sum_{u=-M}^{M} \phi_{i+u}(\xi,\omega_L, \omega_C, \omega_R) q_{i+u}^n \;,
\end{array}
\end{eqnarray}
where $ \omega_L,$ $\omega_C$, $\omega_R$ represents the weights of the WENO reconstruction and they are randomly chosen. Here, $ \phi_{i+u} $ involves the Legendre polynomial and the coefficients of reconstruction polynomials and thus they are different for different orders of accuracy, see \ref{sec:reconstruction} where the shape of these expressions for the third order case are presented. Furthermore, the strategy for obtaining them for other orders of accuracy is also detailed. Then by replacing $q_i^n = A^n \cdot e^{ I \phi \cdot i}$ we obtain
\begin{eqnarray}
\begin{array}{c}
\displaystyle
w_i(\theta,\xi) = A^n  \cdot \sum_{u=-M}^{M} \phi_{i+u}(\xi, \omega_L, \omega_C, \omega_R)  \cdot e^{ (i+u)\theta I} \;,
\end{array}
\end{eqnarray}
that is, once $ \omega_S$ for $S=L,C,R$ are randomly chosen the reconstruction polynomial yields a function of $\theta$ and $\xi$.  Therefore, by collecting all the previous expression, the predictor step takes the form
\begin{eqnarray}
\begin{array}{rl}
q_i(\theta, c, r, \xi, \tau ) &:= 
 w_i(\theta,\xi) 

\\

&
\displaystyle

+ \sum_{k=1}^M  \frac{ ( \tau )^k }{ k!} 

\sum_{l=0}^k 
\sum_{j = 1}^l \left(
\begin{array}{c}

l \\
j

\end{array}
\right) 
\frac{ d^{j} }{dx^{j}} w_i(\theta, \xi) \cdot (-c)^{j} \cdot (r)^{l-j}
\;.
\end{array}
\end{eqnarray}
Then, by inserting this expression into (\ref{eq:amplitude-ADER}) we obtain the corresponding amplitude, which will be denoted by $A^{Explicit}(\theta, c, r)$.  Since $\omega_L$, $\omega_C$ and $\omega_R$ are randomly chosen and due to the factor $A^{-n} \cdot e^{- I \theta i}$ in (\ref{eq:amplitude-ADER}), the amplitude does not depend on $i$.

\subsubsection{Implicit predictor}
Now, let us use the same approach for the implicit predictor, presented in this work. Therefore, the predictor is given by
\begin{eqnarray}
\label{eq:implicit-predictor}
\begin{array}{rl}
q_i(\theta, c, r, \xi, \tau ) &:= 
 w_i(\theta,\xi) 

\\

&
\displaystyle

- \sum_{k=1}^M  \frac{ ( \tau )^k }{ k!} 

\sum_{l=0}^k 
\sum_{j = 1}^l \left(
\begin{array}{c}

l \\
j

\end{array}
\right) 
\partial_x^{(j)} q_i(\theta, c, r, \xi, \tau ) \cdot (-c)^{j} \cdot (r)^{l-j}
\;.
\end{array}
\end{eqnarray}
with 
\begin{eqnarray}
\label{eq:j-x-derivative}
\begin{array}{cc}

\partial_x^{(j+1)} q_i(\theta, c, r, \xi, \tau )

=
&
\partial_x^{(j)} w_i(\theta, c, r, \xi)

 + \tau \cdot \biggl(  
r \cdot \partial_x^{(j)} q_i(\theta, c, r, \xi, \tau )

\\

&
 - c \cdot \partial_x^{(j+1)} q_i(\theta, c, r, \xi, \tau )
\biggr)
\;,

\end{array}
\end{eqnarray}
for $j=1,..,M-1$ and 
\begin{eqnarray}
\label{eq:M-x-derivative}
\begin{array}{cc}

\partial_x^{(M)} q_i(\theta, c, r, \xi, \tau )

=

\partial_x^{(M)} w_i(\theta, c, r, \xi)
 + \tau \cdot   
r \cdot \partial_x^{(M)} q_i(\theta, c, r, \xi, \tau ) 
\;.

\end{array}
\end{eqnarray}

Notice that the predictor in this case is obtained by solving simultaneously (\ref{eq:implicit-predictor})-(\ref{eq:j-x-derivative})-(\ref{eq:M-x-derivative}). Then, by inserting $q_i(\theta, c, r, \xi, \tau )$ into (\ref{eq:amplitude-ADER}) we obtain the corresponding amplitude, which will be denoted by $A^{Implicit}(\theta, c, r)$.  Since $\omega_L$, $\omega_C$ and $\omega_R$ are randomly chosen, the amplitude does not depend on $i$.

Then, the von Neumann analysis for a given order of accuracy $M+1$ is generated as follow:
\begin{itemize}
\item Step 1: Set $c$ and $r$.

\item Step 2: For $j \in \Omega_L \cup \Omega_C \cup \Omega_R$, generate $\tilde{\omega}_L, \tilde{\omega}_C $ and $\tilde{\omega}_R$ as uniformly distributed numbers between $0$ and $1$. Then normalise as $ \omega_S = \frac{ \tilde{\omega_S}}{\omega_L + \omega_C + \omega_R}$, for $S=L,C,R$.

\begin{itemize}

\item Step 2.1:  For $j\in \Omega_L \cup \Omega_C \cup \Omega_R$, compute $ \phi_{i+j}$ as suggested in the appendix \ref{sec:reconstruction}. 

\item Step 2.2: Compute the amplitude $A(\theta, c, r)$ for the corresponding approach, implicit or explicit. 

\end{itemize}

\item Step 3. Repeat Step 2, $NS$ times and obtain the ratio between the number of scenarios in which $|A|<1$ and $NS$. 
     
\end{itemize}

\begin{figure}
\begin{center}
\includegraphics[scale=0.5]{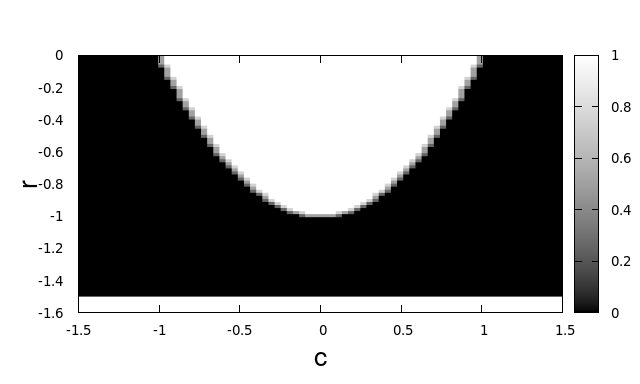}
\end{center}
\caption{Stability for first order scheme, using the numerical flux FORCE. }\label{fig:stab-order-1}
\end{figure}

For numerical implementation, we take $NS = 100$, then we obtain the expected stability region for each order. This will give us meaningful information about the stability region. We implement the FORCE flux scheme and the explicit GRP  whereas for the implicit GRP solver presented here we use the FORCE-$\alpha$ numerical flux. Figure \ref{fig:stab-order-1} shows the result for the first order schemes, the white region is the zone for which the scheme is stable. On the $c$ axis, it is shown the stability for conventional FORCE for the case of purely linear advection equation ($r=0)$, that is, it is reproduced the stability range of $|c|\leq 1$.

\begin{figure}
\begin{center}
\begin{tabular}{cc}
\includegraphics[scale=0.25]{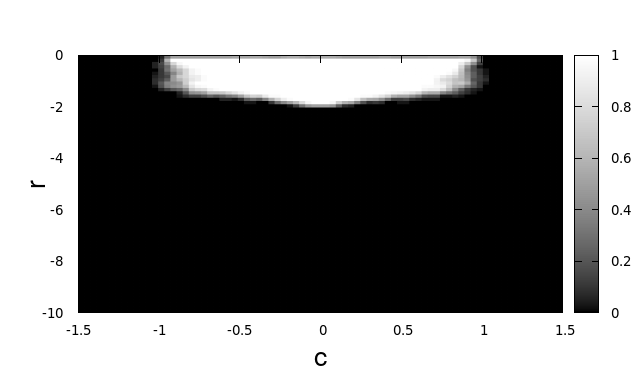}
&
\includegraphics[scale=0.25]{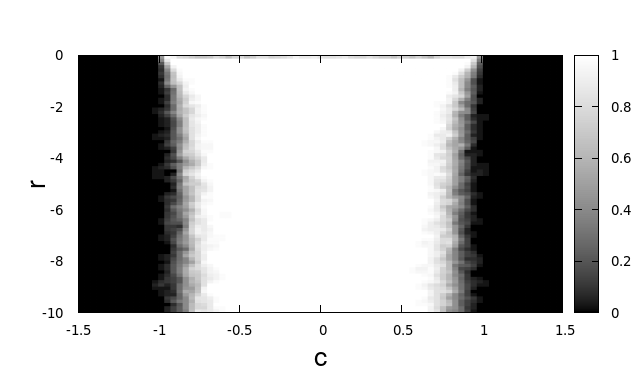} \\
{\small (a) $2nd$ order explicit.} & {\small (b) $2nd$ order implicit.} \\
\\
\includegraphics[scale=0.25]{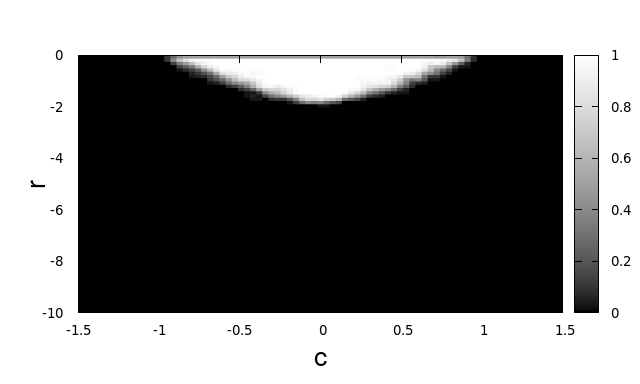}
&
\includegraphics[scale=0.25]{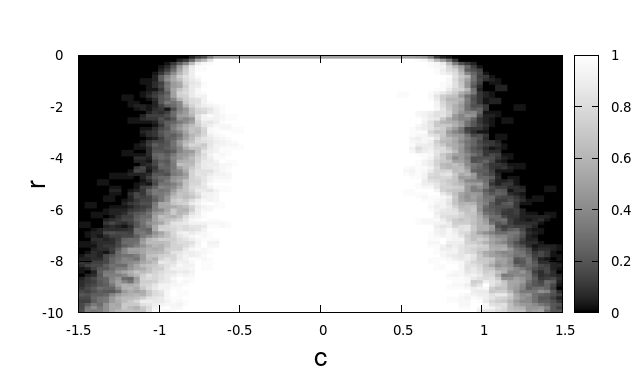}\\
{\small (c) $3rd$ order explicit.} & {\small (d) $3rd$ order implicit.} \\
\\
\includegraphics[scale=0.25]{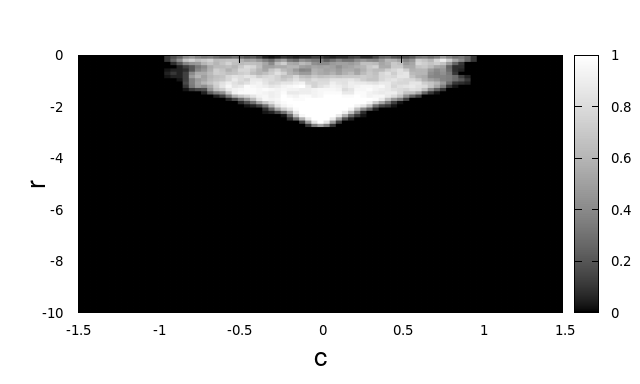}
&
\includegraphics[scale=0.25]{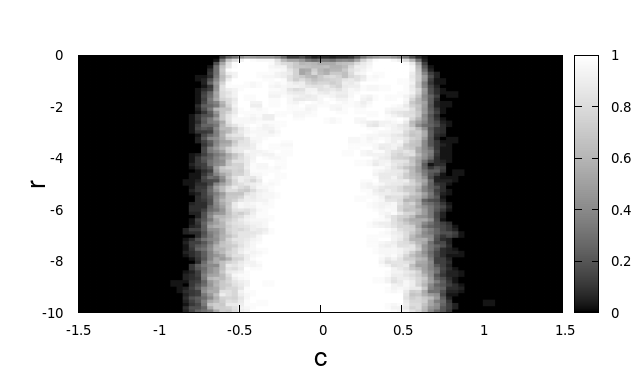}\\
{\small (e) $4th$ order explicit.} & {\small (f) $4th$ order implicit.} \\
\\
\includegraphics[scale=0.25]{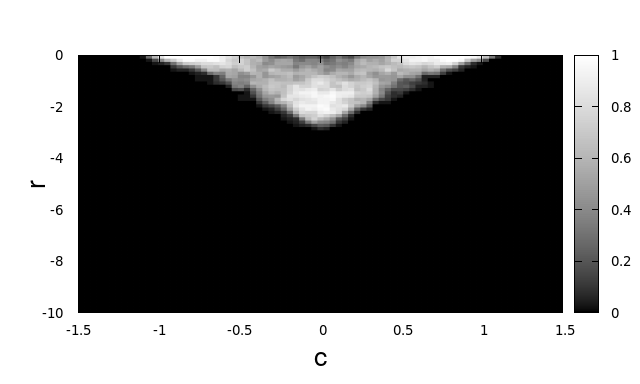}
&
\includegraphics[scale=0.25]{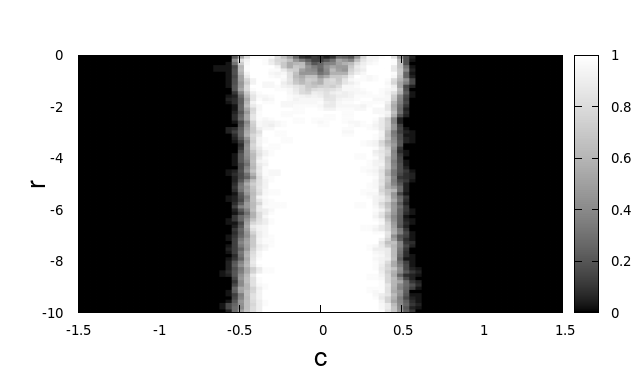}\\
{\small (g) $5th$ order explicit.} & {\small (h) $5th$ order implicit.} \\
\end{tabular}
\end{center}
\caption{Stability of ADER schems. Implicit schemes use FORCE-$\alpha$ with $\alpha = 1$.}\label{fig:stab-high-order}
\end{figure}

Figure \ref{fig:stab-high-order}, shows the stability analysis result for schemes of 2nd, 3rd, 4th and 5th orders of accuracy. On the left part it is shown the result for explicit GRP solver. On the right part it is shown the result for the implicit GRP solver, both figures use the FORCE numerical flux. As can be seen, the range of coefficients $c$ is reduced in comparison with the explicit approach. However, the gaining on the stability for dealing with source terms is evident. The range of $ r $ is limited up to $-10$, however, for large values, that is, in stiff regimes, the scheme still has a range of $c$ which is 
within the range of practical implementations. Since the stability requires a small $C_{CFL}$ coefficient we profit from the new centred, low-dissipation numerical flux in \cite{Toro:2020a}, the FORCE-$\alpha$.  The figure shows the stability range for the case of the implicit scheme, let us call $C_{max}$ the range of $c$ values for which the scheme is stable.  In \cite{Toro:2020a} it is shown that a suitable combination of $ C_cfl$ and $\alpha$ provides a dissipation which is comparable with that of the Godunov method at $C_{max}$.

\begin{figure}
\begin{center}
\begin{tabular}{cc}
\includegraphics[scale=0.25]{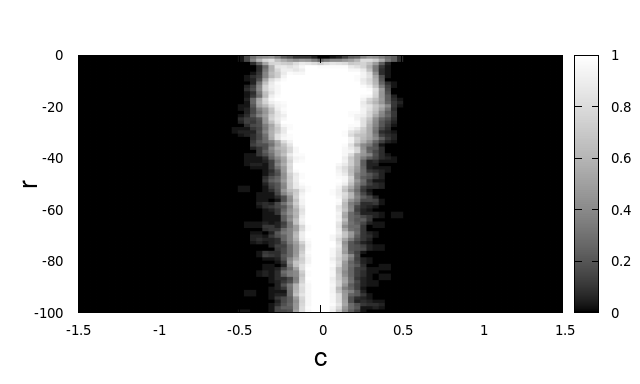}
&
\includegraphics[scale=0.25]{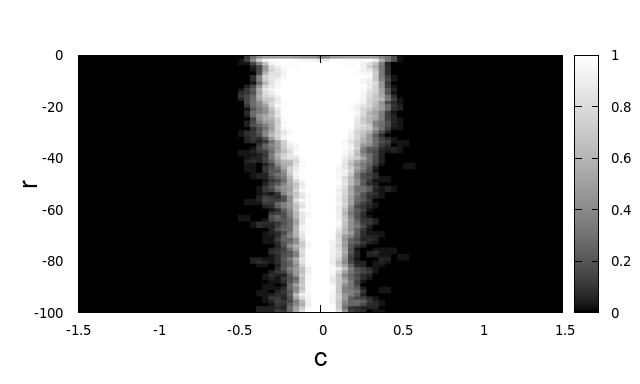} \\
{\small (a) $\alpha = 1$.} & {\small (b)  $\alpha = 1.5$} \\

\includegraphics[scale=0.25]{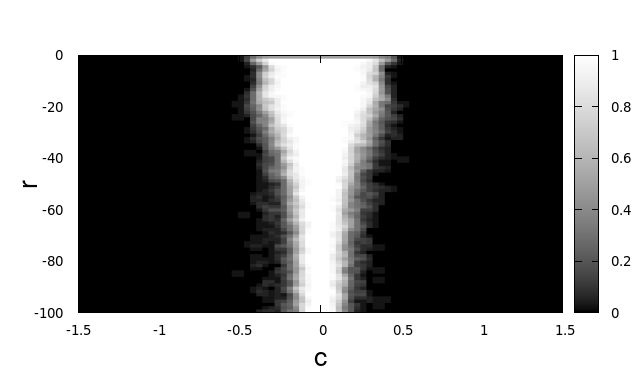}
&
\includegraphics[scale=0.25]{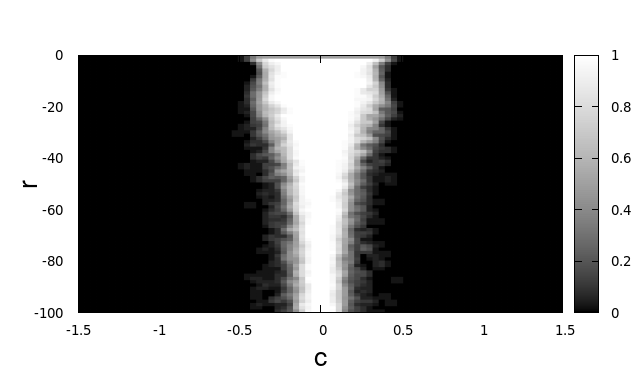} \\
{\small (b) $\alpha = 2$.} & {\small (c)  $\alpha = 5.0$} \\
\\
\includegraphics[scale=0.25]{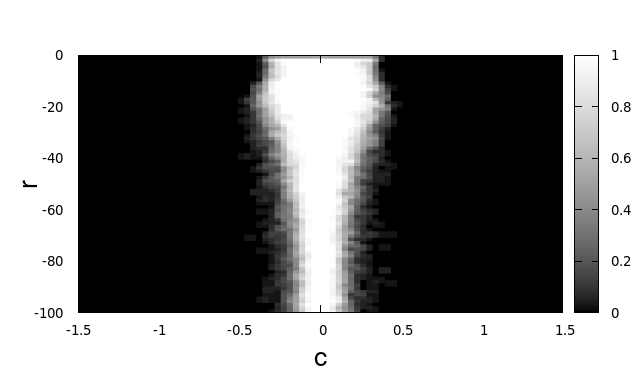}
&
\includegraphics[scale=0.25]{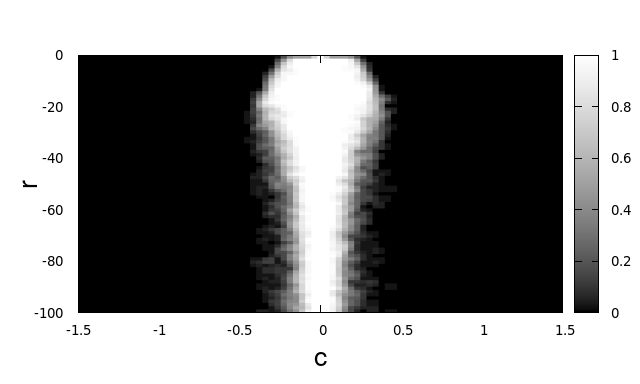}\\
{\small (d)  $\alpha = 10$.} & {\small (e)  $\alpha = 20$.} \\
\\
\includegraphics[scale=0.25]{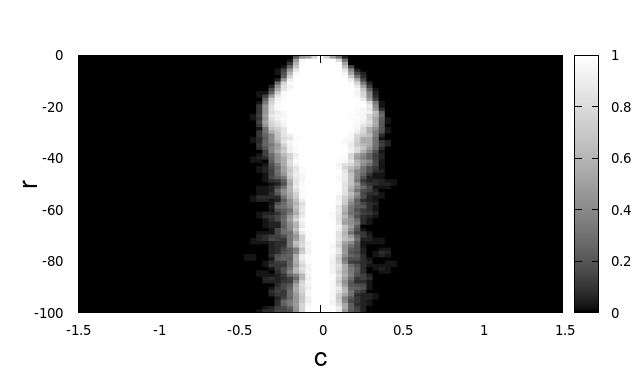}
&
\includegraphics[scale=0.25]{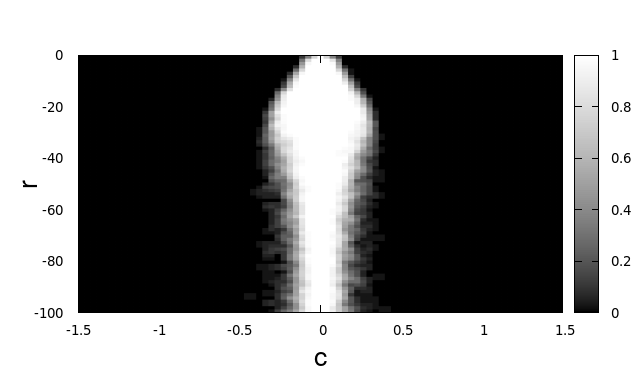}\\
{\small (f)  $\alpha = 50$.} & {\small (g)  $\alpha = 100$.} \\
\end{tabular}
\end{center}
\caption{Stability of the implicit 5$th$ order schems. For different values of $\alpha$ for FORCE-$\alpha$.}\label{fig:stab-several-alpha-5thOrd}
\end{figure}

Figure \ref{fig:stab-several-alpha-5thOrd}, shows the stability for the fifth order scheme for several values of $\alpha$, used in the FORCE-$\alpha$ numerical flux. Notice that the analysis carried out in \cite{Toro:2020a} is not applicable to the scheme in this paper. The scheme in \cite{Toro:2020a} is a first order accurate one and devoted to solve conservation laws. In this case we observe a reduction on the range of $c$ for large values of $\alpha$, particularly for $ \alpha \geq 20 $. 

From the analysis carried out here, we can ensure that $C_{CFL} = 0.1$ belongs to the range of values for which the scheme is stable up to fifth orders of accuracy. Furthermore, this value is comparable with the $C_{CFL}$ range of $
\frac{1}{2N+1}$ for the Discontinuous Galerkin approach. Despite, the range can be considered small, as shown in the next section the compromise between accuracy and efficiency is good enough to make this scheme a feasible high-order method.

\subsection{An efficiency analysis associated to the present scheme}\label{sec:operations}

In this section, a comparison of the order of operations between the present approach and that in \cite{Toro:2015a} for obtaining the predictor is carried out.  

Both approaches have in common the following steps to get the predictor.

\begin{itemize}
\item Composition of a non-linear algebraic system $\mathcal{H}$ which involves implicit Taylor series expansions to compute the state and the evolution of the spatial derivatives as well.

\item The composition of the Jacobian matrix $\nabla \mathcal{H}$ of the non-linear algebraic system $\mathcal{H}$.

\item A descent step for fixed point iteration procedures is implemented.  That is, the sought state $\mathbf{Q}$ is computed iteratively as $\mathbf{Q}^{r+1} = \mathbf{Q}^r-\delta$, where $\delta$ solves the linear system $\nabla \mathcal{H} \delta = \mathcal{H}$.  This involves an inversion matrix procedure which is carried out once $\mathcal{H}$ and $\nabla \mathcal{H}$ have been constructed. 

\end{itemize}

To differentiate both procedures, we are going to denote by $O(\mathcal{H}^{MT})$, $O(\nabla \mathcal{H}^{MT})$ and $O(\delta^{MT})$ the orders of the number of operations to get these operators in the case of the $MT$ solver in \cite{Toro:2015a}. Similarly, we are going to denote by $O(\mathcal{H}^{MTS})$, $O(\nabla \mathcal{H}^{MTS})$ and $O(\delta^{MTS})$ the 
orders of the operations to get these operators in the case of the present approach.

Let us assume that the Cauchy-Kowalewskaya procedure can be written as
\begin{eqnarray}
\begin{array}{c}
\partial_t^{(k)} \mathbf{Q} = \mathbf{G}^{k}=  \sum_{l=0}^{k} \mathbf{C}_{l,k} \cdot \partial_x^{(l)} \mathbf{Q}\;, 
\end{array}
\end{eqnarray}
where $\mathbf{C}_{l,k}$ is a $m\times m$ coefficient matrix which may depend on the data itself. Furthermore, base on the linear system equations we can assume that the number of operations to build $\mathbf{C}_l^{k} \cdot \partial_x^{(l)}$  is $ O( \mathbf{C}_l^{k} \cdot \partial_x^{(l)} ) =O( k m^{k+1})$. Therefore, the Taylor series for computing the data $\mathbf{Q}$, yields
\begin{eqnarray}
\begin{array}{c}
O(\mathbf{Q}) = \sum_{k=1}^{M} O(\mathbf{G}^{k}) = O(\sum_{k=1}^M  k m^{k+1} ) \;.
\end{array}
\end{eqnarray}
Since 
\begin{eqnarray}
\label{eq:summation-kmk}
\begin{array}{c}
 \sum_{k=1}^M k m^{k-1} = \frac{d}{dm}( \sum_{k=1}^M m^{k}) = \frac{d}{dm}(\frac{m^{M}-1}{m-1})= \frac{1 + m^{M}(M-1) - Mm^{m-1}  }{ (m-1)^2}
 \;,
\end{array}
\end{eqnarray}
one obtains
\begin{eqnarray}
\begin{array}{c}
O(\mathbf{Q}) =  O(M m^{M} ) \;,
\end{array}
\end{eqnarray}

On the other hand, in \cite{Toro:2015a} the evolution of the spatial derivatives is obtained by mean of the implicit Taylor series expansions
\begin{eqnarray}
\begin{array}{c}
\partial_x^{(l)} \mathbf{Q}
=
\mathbf{w}_x^{(l) }- \sum_{k=1}^{M-l} \frac{(-\tau)^k}{k!} \partial_t^{(k)} ( \partial_x^{(l)} \mathbf{Q} )
\;.
\end{array}
\end{eqnarray}
We can express this expansion as
\begin{eqnarray}
\begin{array}{c}
\partial_x^{(l)} \mathbf{Q}
=
\mathbf{w}_x^{(l) }- \sum_{k=1}^{M-l} \frac{(-\tau)^k}{k!} \partial_t^{(k)} ( \partial_x^{(l)} \mathbf{Q} )
=
\mathbf{w}_x^{(l) } + \sum_{k=0}^{M} \mathbf{D}_{l,k} \partial_x^{(k)}\mathbf{Q}
\end{array}
\end{eqnarray}
where $\mathbf{D}_{l,k}$ is a $m\times m$ matrix, they can depend on the data it self. We assume that each one of the expressions $\mathbf{D}_{l,k} \partial_x^{(k)}\mathbf{Q}$ involves $k m^{k+1}$ operations. Therefore
\begin{eqnarray}
\begin{array}{c}
O(\partial_x^{(l)} \mathbf{Q}) = O(\sum_{k=0}^{M-l} \mathbf{D}_{l,k} \partial_x^{(k)}\mathbf{Q}) = O( \sum_{k=1}^{M-l} k m^{k+1})

=O( (M-l)m^{M-l}) \;.
\end{array}
\end{eqnarray}
Then, to compose $\mathcal{H}^{MT}$ one needs
\begin{eqnarray}
\begin{array}{c}
O( \mathcal{H}^{MT}) = O( \mathbf{Q}) +  \sum_{l=1}^{M} O( \partial_x^{(l) \mathbf{Q}} ) 
=
O(M m^m) + O( \sum_{l=1}^{M} (M-l)m^{M-l} )
\end{array}
\end{eqnarray}
operations. 
From (\ref{eq:summation-kmk}), $O( \sum_{l=1}^{M} (M-l)m^{M-l}) = O( (M-1)\cdot m^{M-1} $ and thus 
\begin{eqnarray}
\begin{array}{c}
O( \mathcal{H}^{MT}) = O( M m^{M} )\;. 
\end{array}
\end{eqnarray}
To build the Jacobian $ \nabla \mathcal{H}$ we note that 
\begin{eqnarray}
\begin{array}{c}
\displaystyle
(\nabla H)_{0,s} = \sum_{k=0}^{M} \frac{ \partial (\mathbf{C}_{l,k} \partial_x^{(k)} \mathbf{Q} ) }{ \partial_{ \partial_x^{(s)} \mathbf{Q}}}
\end{array}
\end{eqnarray}
and 
\begin{eqnarray}
\begin{array}{c}
\displaystyle
(\nabla H)_{l,s} = \sum_{k=0}^{M} \frac{ \partial (\mathbf{D}_{l,k} \partial_x^{(k)} \mathbf{Q} ) }{ \partial_{ \partial_x^{(s)} \mathbf{Q}}} \;,
\end{array}
\end{eqnarray}
for $l=1,...,M$.
Then, let us assume that 
\begin{eqnarray}
\begin{array}{c}
O(\frac{ \partial (\mathbf{C}_{l,k} \partial_x^{(k)} \mathbf{Q} ) }{ \partial_{ \partial_x^{(s)} \mathbf{Q}}})
=  k m^{k+1}
\;,
O(\frac{ \partial (\mathbf{D}_{l,k} \partial_x^{(k)} \mathbf{Q} ) }{ \partial_{ \partial_x^{(s)} \mathbf{Q}}})
= O(k m^{M} )\;
\end{array}
\end{eqnarray}
and by noting that $\nabla \mathcal{H}$ is a $(M+1) \times (M+1)$ block matrix. 

\begin{eqnarray}
\begin{array}{c}

\displaystyle
O((\nabla H)_{l,s} ) = O(  \sum_{k=0}^{M} k m^M ) = O( m^M  \cdot \frac{M(M+1)}{2}) = O(M^2 \cdot m^M)\;.  
\end{array}
\end{eqnarray}

Then the number of operations to form $\nabla \mathcal{H}^{MT}$ can be assumed to be of the order 
\begin{eqnarray}
\begin{array}{c}
O(\nabla \mathcal{H}^{MT}) =  O( M^4 \cdot m^{M}) \;.
\end{array}
\end{eqnarray}
Notice that once the Jacobian matrix $ \nabla \mathcal{H}$ is computed the inversion matrix involves $m^3 M^3$ operations which is neglected by comparison with the number of operations for $\mathcal{H}$ and $\nabla \mathcal{H}$. Let us define the number of operations for computing the predictor by using \cite{Toro:2015a} as 
\begin{eqnarray}
\begin{array}{c}
O(MT) = O(\mathcal{H}^{MT}) + O(\nabla \mathcal{H}^{MT})=  O(m^{M}(M^4 + M)) \;.
\end{array}
\end{eqnarray}
Now, let us compute the order of operations involved in the present approach.  In this case $O(\mathcal{H}^{MTS})$ involves the computation of the state, which by following the same procedure above yields $ O(\mathbf{Q}) = O(M m^m )$.  Since all the derivatives involves a second order evolutions $O(\partial_x^{(l)} \mathbf{Q}) = O(2 m^2)$. In the present approach, the non-linear algebraic systems is only for the implicit Taylor series approximation for $\mathbf{Q} $, then $ O( \nabla \mathcal{H} ) = O(M m^M)$.  Therefore the order of operations to obtain the predictor in the present approach is assumed to be
\begin{eqnarray}
\begin{array}{c}
O(MTS) = O(M m^{M})\;.
\end{array}
\end{eqnarray}

On the other hand, it is expected to have a penalization on the efficiency of the current approach due to the small CFL ranges ($C_{cfl} = 0.1$). The factor in which performance is reduced compared to the approach in \cite{Toro:2015a} which uses a $CFL$ coefficient of $C_{cfl}^{MT} = 0.9$ should be $\Delta_{CFL} := \frac{ C_{cfl}^{MT}}{C_{cfl}^{MT}} = 9\;.$

Despite the current scheme does not use numerical fluxes, but jump operators which are expected to behave as a conservative scheme if the PDE admits a conservation form, we can assume that both procedures in terms of performances only have a difference in the way in which predictors have been computed. Then we can approximate the gaining factor in terms of the performance of the present approach against that  in \cite{Toro:2015a} by means of
\begin{eqnarray}
\begin{array}{c}
\frac{O(MTS)}{O(MT)} \cdot \Delta_{CFL} =  O( \frac{9}{ M^3+1})\;.
\end{array}
\end{eqnarray}
So, we define the Gaining coefficient $O(G_M)$ give by
\begin{eqnarray}
\label{eq:gaining:coefficient}
\begin{array}{c}
G_M = \frac{9}{ M^3+1}\;.
\end{array}
\end{eqnarray}
here $\frac{O(MTS)}{O(MT)}  $ represents the gaining factor due to the difference on the predictor step.  Notice that this coefficient does not depend on the number of variables $m$. Therefore, we expect the present scheme to be less efficient than the scheme in \cite{Toro:2015a} when $G_{M} > 1$. Remember that $M+1$ corresponds to the expected order of accuracy of the scheme. So, the gaining factor predicts that the present scheme becomes the most efficient (related to that in \cite{Toro:2015a} and under the condition of $C_{CFL}=0.1$) from the third order onwards. However, the second order is exactly the same for both approaches and the present one works well with $C_{CFL}^{MTS} = 0.9$, the associated penalization disappears and then the scheme is as efficient as the second order scheme in \cite{Toro:2015a}.

As the number of operations is directly related to the CPU time, we can empirically measure the gaining in performance of the scheme, just by taking the quotients between the CPU time for the scheme in \cite{Toro:2015a} and the present one. Let us make the comparison for a linear system case in which $m = 2$. This is because the orders of operations are consistent with that used to derive the coefficient (\ref{eq:gaining:coefficient}).
\begin{table}
\begin{center}
Theoretical order : 2, $M = 1$  
\\
\begin{tabular}{cccccccc} 
\\ 
\hline
\hline
          & 128 cells &    64  cells & 32 cells & 16 cells \\
\hline          
CPU - MTS &  $ 2.84 $ & $ 0.78 $     & $ 0.25 $  &  $ 0.24 $ \\
CPU - MT  &  $ 0.89 $ & $ 0.23 $     & $ 0.05 $  &  $ 0.17 $ \\
\hline
\rowcolor{Gray} CPU ratio & 3.19      & 3.39    &  5      &   1.14 \\
\hline
Expected   &  gaining          & $ G_1 = 4.5 $       &           &     \\
\hline
\\
\end{tabular} 
\\
Theoretical order : 3, $M = 2$  
\\
\begin{tabular}{cccccccc} 
\\
\hline
\hline 

            & 128 cells & 64   cells & 32 cells & 16 cells \\
\hline          
CPU - MTS   & $ 9.36  $ & $ 2.51 $ & $ 0.73 $ & $ 0.17 $  \\
CPU - MT    & $10.06  $ & $ 2.74 $ & $ 0.77 $ & $ 0.21 $ \\
\hline
\rowcolor{Gray} CPU ratio   & $ 0.93 $  & $ 0.92 $ & $ 0.95 $ &  $ 0.81$ \\
\hline
Expected   &  gaining          & $ G_2 = 1$       &           &     \\
\hline
\\
\end{tabular} 
\\
Theoretical order : 4, $M = 3$  
\\
\begin{tabular}{cccccccc} 
\\
\hline 

            & 128 cells & 64 cells & 32 cells & 16 cells \\
\hline          
CPU - MTS   & 27.32     & 6.99     &  1.80     &   0.50  \\
CPU - MT    & 108.56    & 27.55    & 7.46     &   1.72  \\
\hline
\rowcolor{Gray} CPU ratio   & 0.25      & 0.25    &  0.24      &   0.29 \\
\hline
Expected   &  gaining          & $ G_3 = 0.31 $       &           &     \\
\hline
\\
\end{tabular} 
\\
Theoretical order : 5, $M = 4$  
\\
\begin{tabular}{cccccccc} 
\\
\hline
\hline 
            & 128 cells & 64  cells & 32 cells & 16 cells \\
\hline          
CPU - MTS   &  89.01  &  22.51 &  5.77 & 1.55  \\
CPU - MT    & 671.68  & 169.84 & 42.88 & 10.71 \\
\hline
\rowcolor{Gray} 
CPU ratio   &   0.13      & 0.13    &  0.13      &   0.14 \\
\hline
Expected    &  gaining          & $ G_4 = 0.14 $       &           &     \\
\hline
\end{tabular}

\end{center}
\caption{Linear system. Output time
  $t_{out} = 1$ with, $C_{cfl}^{MT}= 0.9$ and $C_{cfl}^{MTS}= 0.1$, $\alpha = 1.9$, $\beta = -1$, $\lambda =  1 $.}\label{Table-GM-LinearSystem}
\end{table}

Table \ref{Table-GM-LinearSystem} shows the CPU time comparison between the present scheme (MTS) and that in \cite{Toro:2015a} (MT). For the orders 2nd, 3rd, 4th and 5th, the columns show the CPU times for $128$, $64$, $32$ and $16$ computational cells.
The row {\it CPU ratio} shows the quotient between both CPU times for each number of cells. This coefficient measures empirically the gaining on efficiency. Furthermore, the values are in the range of the expected values predicted by $G_M$.
As expected, since the second order has been used with $C_{CFL} = 0.1$, this is less efficient that the counterpart in  \cite{Toro:2015a}.

The gaining factor is expected to improve for the non-linear case. However, just in the scalar linear advection case, the Jacobian $\nabla \mathcal{H}$ of the non-linear equation for \cite{Toro:2015a} is a diagonal matrix and  $O(\nabla \mathcal{H}) = O(M^2) $. Therefore, the gaining factor becomes 
\begin{eqnarray}
\begin{array}{c}
G_M =O (\frac{9}{M}) \;,
\end{array}
\end{eqnarray}
this analytically demonstrates that that is the only case in which the preset approach is less efficient than the scheme in \cite{Toro:2015a}.

\section{Numerical Results}\label{sec:num-results}
In this section, a set of test problems is considered aimed at proving the applicability of the present scheme.
In all the implementations, the time step is computed as
\begin{eqnarray}
\begin{array}{c}
\Delta t = C_{cfl} \frac{ \Delta x}{ \lambda_{max} }\;,
\end{array}
\end{eqnarray}
where $ \lambda_{max} = \max_{i=1,...,N}\{ max_{j = 1,...,m}\{ |\lambda_j( \mathbf{Q}_i^n) \} \}$, with $ \lambda_j(\mathbf{Q})$  is the $j-th$ eigenvector of the Jacobian matrix $\mathbf{A}(\mathbf{Q})$.   In this section, we

\subsection{The LeVeque and Yee test}

\begin{figure}
\begin{center}
\includegraphics[scale=0.5]{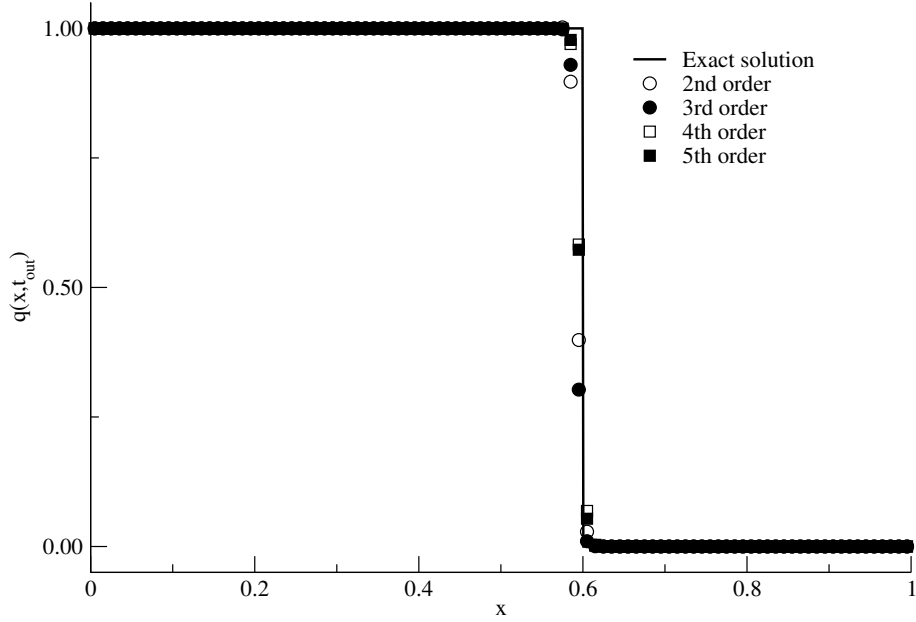}
\end{center}
\caption{Leveque and Yee test. We have used $100$ cells, $C_{cfl}=0.1$, $\alpha=2.4$, $t_{out} = 0.3$, $\beta  = -1000$.}\label{fig:LevequeAndYee}
\end{figure}

Here, we apply the present schemes to the well-known and challenging scalar test
problem proposed by LeVeque and Yee \cite{LeVeque:1990a}, given by
\begin{eqnarray}
\begin{array}{c}
\partial_t q(x,t) + \partial_x q(x,t) = \beta q(x,t)(q(x,t) - 1)(q(x,t) - \frac{1}{2} ) \;.
\end{array}
\end{eqnarray}
We solve this PDE on the computational domain $[0,1]$ with transmissive boundary conditions and the initial condition given by
\begin{eqnarray}
\begin{array}{c}
q(x,0) = \left\{

\begin{array}{cc}
1 \;, x < 0.3 \;, \\
0 \;, x > 0.3 \;. \\
\end{array}
\right.
\end{array}
\end{eqnarray}
The solution on the characteristic curves satisfies de ordinary differential equation  $ \frac{d(x(t), t)}{dt} = \beta q(x(t),t)(q(x(t),t) - 1)(q(x(t),t) - \frac{1}{2} )  $, which has two stable solutions $q \equiv 0$ and $q \equiv 1$ and one unstable solution in $q \equiv \frac{1}{2}$ where any solution trays to away from this. Similarly, any solution associated with characteristic curves necessarily must converge to one of the two stable solutions. On the other hand, a numerical scheme that is not able to solve stiff source terms may introduce an excessive numerical diffusion and so the numerical solution, following characteristic curves, converges to the wrong stable solution. This penalizes the right propagation. Figure \ref{fig:LevequeAndYee} shows the comparison between the exact solution and the numerical approximations provided by the present scheme of second, third, fourth and fifth orders of accuracy. The figure shows a good agreement for $\beta = - 1000$ at $t_{out} = 0.3$, which correspond to the stiff regime. We have used $ C_{cfl}= 0.1$ and $100$ cells. This test illustrates the ability of the present scheme for solving hyperbolic balance laws with stiff source terms.
%

\subsection{A linear system of hyperbolic balance laws}
Here, we consider the linear system in \cite{Montecinos:2014c}, given by
\begin{eqnarray}
\begin{array}{c}
\partial_t \mathbf{Q}(x,t) + \mathbf{A} \partial_x ( \mathbf{Q} (x,t)) = \mathbf{B} \mathbf{Q} (x,t) \;, x \in [0,1]\;,\\
\mathbf{Q}(x,0) = 
\left[ 
\begin{array}{c}

\sin( 2 \pi x ) \\
\cos( 2 \pi x ) 

\end{array}
\right]
\;,
\end{array}
\end{eqnarray}
where 
\begin{eqnarray}
\begin{array}{c}

\mathbf{A} = \left[
\begin{array}{cc}
0       & \lambda \\
\lambda & 0 \\
\end{array}
\right]\;,

\mathbf{B} = \left[
\begin{array}{cc}
\beta   & 0 \\
0 & \beta \\
\end{array}
\right]
\;.
\end{array}
\end{eqnarray}

The problem is endowed with periodic boundary conditions.  
This system has the exact solution 
\begin{eqnarray}
\begin{array}{c}
\displaystyle \mathbf{Q}^{e}(x,t) = 
\frac{  e^{\beta t} }{2} 
\left[
\begin{array}{c}
\Phi(x,t ) + \Psi(x,t )  \\
\Phi(x,t ) - \Psi(x,t )  \\
    
\end{array}
\right] \;,
\end{array}
\end{eqnarray}
where
\begin{eqnarray}
\begin{array}{c}

\Phi(x,t ) = \sin( 2\pi (x - \lambda t) ) + \cos( 2\pi (x - \lambda t) )  \;, \\

\Psi(x,t) = \sin( 2\pi (x + \lambda t) ) - \cos( 2\pi (x + \lambda t) )\;.
\end{array}
\end{eqnarray}
Here we consider $\lambda =1$, $\alpha = 1.9$, $\beta = -1$ and $C_{cfl}=0.1$. This is a simple test aimed at evaluating the accuracy of the present scheme. As can be seen in Table \ref{Table-LinearSystem}, the expected theoretical orders of accuracy are achieved.  Notice that, the second order is super optimal, this is because for a such order the scheme coincides with that in \cite{Toro:2015a} and the accuracy limit is $C_{cfl} = 0.9$ and thus a reduction of this coefficient to $C_{cfl}=0.1$ penalizes the efficiency but increases the accuracy.
\begin{table}
\begin{center}
Theoretical order : 2 \\
\begin{tabular}{cccccccc} 
\\
\hline
\hline 
Mesh  & $L_\infty$ - err & $L_\infty$- ord  & $L_1$ - err & $L_1$ - ord & $L_2$ - err & $L_2$ - ord & CPU  \\  
\hline
    
    16  &  $   5.38e-0 2  $  &   2.36  &  $   3.59e-0 2  $  &   2.43  &  $   3.84e-0 2  $  &   2.45  &   0.2379  \\
    32  &  $   6.34e-0 3  $  &   3.09  &  $   3.77e-0 3  $  &   3.25  &  $   3.93e-0 3  $  &   3.29  &   0.2514  \\
    64  &  $   6.52e-0 4  $  &   3.28  &  $   3.08e-0 4  $  &   3.61  &  $   3.50e-0 4  $  &   3.49  &   0.7839  \\
   128  &  $   1.15e-0 4  $  &   2.50  &  $   1.58e-0 5  $  &   4.29  &  $   3.41e-0 5  $  &   3.36  &   2.8396  \\    

\hline
 \\
 \end{tabular}  
\\
Theoretical order : 3 \\
\begin{tabular}{cccccccc}  
\\
\hline 
Mesh  & $L_\infty$ - err & $L_\infty$- ord  & $L_1$ - err & $L_1$ - ord & $L_2$ - err & $L_2$ - ord & CPU  \\  
\hline
    
    16  &  $   2.88e-0 2  $  &   2.45  &  $   1.88e-0 2  $  &   2.56  &  $   2.08e-0 2  $  &   2.53  &  0.1696  \\
    32  &  $   3.81e-0 3  $  &   2.92  &  $   2.44e-0 3  $  &   2.95  &  $   2.71e-0 3  $  &   2.94  &  0.7285  \\
    64  &  $   4.81e-0 4  $  &   2.99  &  $   3.07e-0 4  $  &   2.99  &  $   3.40e-0 4  $  &   2.99  &  2.5118  \\
   128  &  $   6.02e-0 5  $  &   3.00  &  $   3.83e-0 5  $  &   3.00  &  $   4.26e-0 5  $  &   3.00  &  9.3517  \\
   
\hline
 \\
 \end{tabular}  
\\
Theoretical order : 4 \\
\begin{tabular}{cccccccc}  
\\
\hline 
Mesh  & $L_\infty$ - err & $L_\infty$- ord  & $L_1$ - err & $L_1$ - ord & $L_2$ - err & $L_2$ - ord & CPU  \\  
\hline
    16  &  2.92  &$  1.02e-0 3$&  3.15  &$  5.42e-0 4$&  3.10  &$  6.14e-0 4$  &0.580\\
    32  &  3.64  &$  8.15e-0 5$&  3.67  &$  4.26e-0 5$&  3.67  &$  4.82e-0 5$  &1.8017\\
    64  &  3.87  &$  5.56e-0 6$&  3.87  &$  2.92e-0 6$&  3.87  &$  3.30e-0 6$  &6.9937\\
   128  &  3.95  &$  3.59e-0 7$&  3.95  &$  1.89e-0 7$&  3.95  &$  2.14e-0 7$  &27.3153\\

 \hline
 \\
 \end{tabular}  
\\
Theoretical order : 5 \\
\begin{tabular}{cccccccc}  
\\
\hline 
Mesh  & $L_\infty$ - err & $L_\infty$- ord  & $L_1$ - err & $L_1$ - ord & $L_2$ - err & $L_2$ - ord & CPU  \\  
\hline

    16  &  $   8.99e-0 4  $  &   4.75  &  $   5.87e-0 4  $  &   4.86  &  $   6.48e-0 4  $  &   4.84  &   1.5523  \\
    32  &  $   2.66e-0 5  $  &   5.08  &  $   1.71e-0 5  $  &   5.10  &  $   1.89e-0 5  $  &   5.10  &   5.7690  \\
    64  &  $   5.51e-0 7  $  &   5.60  &  $   3.51e-0 7  $  &   5.60  &  $   3.90e-0 7  $  &   5.60  &  22.5083  \\
   128  &  $   1.86e-0 8  $  &   4.89  &  $   1.18e-0 8  $  &   4.89  &  $   1.32e-0 8  $  &   4.89  &  89.0071  \\

 \hline
\end{tabular} 
\end{center}
\caption{Linear system. Output time
  $t_{out} = 1$ with  $C_{cfl}= 0.1$, $\alpha = 1.9$, $\beta = -1$, $\lambda =  1 $.}\label{Table-LinearSystem}
\end{table}
%

\subsection{A non-conservative system of partial differential equations}
In this section, we solve the following hyperbolic system
\begin{eqnarray}
\begin{array}{c}
\partial_t u + \lambda \partial_x u  + u\partial_x v = 2 \pi u(u-1) \;,  \\
\partial_t v + \lambda \partial_x v  + \partial_x u = -2 \pi (v-1) \;,  \\
\end{array}
\end{eqnarray}
with the initial condition $u(x,0) = 1 + \varepsilon \cos(2 \pi x)$ and $v(x,0) = 1 + \varepsilon \sin(2 \pi x)$.
Notice that this system can be written in the matrix form (\ref{eq:gov-equation-non-cons}), with $\mathbf{Q} = [u,v]^T$, with
\begin{eqnarray}
\begin{array}{c}
\mathbf{A}( \mathbf{Q}) =
\left[
\begin{array}{cc}
\lambda & u \\
1       & \lambda
\end{array}
\right]\;,

\mathbf{S}( \mathbf{Q}) =
\left[
\begin{array}{c}
  2 \pi u (u-1) \\
- 2 \pi (v-1)
\end{array}
\right]\;,

\end{array}
\end{eqnarray}
the eigenvalues of $ \mathbf{A}$ are giving by $ \lambda_1 = \lambda - \sqrt{u} $ and $ \lambda_2 = \lambda + \sqrt{u} $. For the given initial condition  and periodic boundary conditions on the interval $[0,1]$, we have that
\begin{eqnarray}
\begin{array}{c}
u(x,t) = 1 + \varepsilon \cos( 2 \pi (x- \lambda t))\;, \\
v(x,t) = 1 + \varepsilon \sin( 2 \pi (x- \lambda t))\;, \\
\end{array}
\end{eqnarray}
is the exact solution of this problem.  Notice that the system cannot be written in a conservative form and thus the universal scheme is a suitable one, it is because no modification has to be done on the present scheme to solve this test.


The numerical solutions are computed up to $t_{out} = 1$ by using $\lambda = 1$, $\varepsilon= 0.02$ and $16$ cells. Table \ref{Table-non-cons-u} shows the error for the variable $u$. As can bee seen, the expected theoretical order of accuracy are achieved.




%
\begin{table}
\begin{center}
Theoretical order : 2 \\
\begin{tabular}{cccccccc} 
\\
\hline
\hline 
Mesh  & $L_\infty$ - err & $L_\infty$- ord  & $L_1$ - err & $L_1$ - ord & $L_2$ - err & $L_2$ - ord & CPU  \\  
\hline

    16  &  $   8.42e-0 3  $  &   1.40  &  $   4.21e-0 3  $  &   1.51  &  $   4.89e-0 3  $  &   1.51  &  0.34  \\
    32  &  $   1.93e-0 3  $  &   2.13  &  $   9.29e-0 4  $  &   2.18  &  $   1.09e-0 3  $  &   2.17  &  1.19  \\
    64  &  $   4.19e-0 4  $  &   2.20  &  $   1.97e-0 4  $  &   2.23  &  $   2.35e-0 4  $  &   2.21  &  4.62  \\
   128  &  $   9.62e-0 5  $  &   2.12  &  $   4.47e-0 5  $  &   2.14  &  $   5.38e-0 5  $  &   2.12  & 13.99  \\
        																															
\hline
\\
\end{tabular}  
\\
Theoretical order : 3 \\
\begin{tabular}{cccccccc}  
\\
\hline 
Mesh  & $L_\infty$ - err & $L_\infty$- ord  & $L_1$ - err & $L_1$ - ord & $L_2$ - err & $L_2$ - ord & CPU  \\  
\hline

    16  &  $   2.40e-0 3  $  &   2.26  &  $   1.42e-0 3  $  &   2.20  &  $   1.58e-0 3  $  &   2.22  &  0.55  \\
    32  &  $   3.45e-0 4  $  &   2.80  &  $   2.09e-0 4  $  &   2.76  &  $   2.33e-0 4  $  &   2.76  &  3.03  \\
    64  &  $   4.52e-0 5  $  &   2.93  &  $   2.77e-0 5  $  &   2.92  &  $   3.08e-0 5  $  &   2.92  & 11.09  \\
   128  &  $   5.75e-0 6  $  &   2.97  &  $   3.55e-0 6  $  &   2.97  &  $   3.94e-0 6  $  &   2.97  & 47.08  \\

 \hline
\\
\end{tabular}  
\\
Theoretical order : 4 \\
\begin{tabular}{cccccccc}  
\\
\hline 
Mesh  & $L_\infty$ - err & $L_\infty$- ord  & $L_1$ - err & $L_1$ - ord & $L_2$ - err & $L_2$ - ord & CPU  \\  
\hline

    16  &  $   9.44e-0 4  $  &   3.50  &  $   5.84e-0 4  $  &   3.65  &  $   6.46e-0 4  $  &   3.62  &   2.45 \\
    32  &  $   4.89e-0 5  $  &   4.27  &  $   2.84e-0 5  $  &   4.36  &  $   3.17e-0 5  $  &   4.35  &   8.94 \\
    64  &  $   2.52e-0 6  $  &   4.28  &  $   1.42e-0 6  $  &   4.33  &  $   1.59e-0 6  $  &   4.32  &  37.31 \\
   128  &  $   1.38e-0 7  $  &   4.19  &  $   7.62e-0 8  $  &   4.22  &  $   8.58e-0 8  $  &   4.21  & 168.41 \\

 \hline
\\
\end{tabular}  
\\
Theoretical order : 5 \\
\begin{tabular}{cccccccc}  
\\
\hline 
Mesh  & $L_\infty$ - err & $L_\infty$- ord  & $L_1$ - err & $L_1$ - ord & $L_2$ - err & $L_2$ - ord & CPU  \\  
\hline

    16  &  $   7.61e-0 5  $  &   4.52  &  $   4.47e-0 5  $  &   4.42  &  $   4.98e-0 5  $  &   4.44  &   6.14   \\
    32  &  $   2.62e-0 6  $  &   4.86  &  $   1.58e-0 6  $  &   4.83  &  $   1.76e-0 6  $  &   4.83  &  26.78  \\
    64  &  $   8.42e-0 8  $  &   4.96  &  $   4.87e-0 8  $  &   5.02  &  $   5.44e-0 8  $  &   5.01  &  92.01  \\
   128  &  $   3.82e-0 9  $  &   4.46  &  $   1.49e-0 9  $  &   5.03  &  $   1.82e-0 9  $  &   4.90  & 229.42  \\
   
 \hline
\end{tabular} 
\end{center}
\caption{Non-conservative system, variable $u$. Output time
  $t_{out} = 1.0 s$ with  $C_{cfl}= 0.1,$ $\lambda = 1$ and $\alpha = 2.2$. }\label{Table-non-cons-u}
\end{table}
%

\subsection{The Euler equations}
Now let us consider the Euler equations, given by
\begin{eqnarray}
\label{eq:1-euler}
\begin{array}{ccc}
\mathbf{Q} =
\left[
\begin{array}{c}
\rho \\
\rho u \\
E
\end{array}
\right]\;,
&
\mathbf{F}(\mathbf{Q}) =
\left[
\begin{array}{c}
\rho u \\
\rho u^2 + p \\
u(E+p)
\end{array}
\right]\;,
\end{array}
\end{eqnarray}
where the pressure $p$ is related with the conserved variables through the equation 
\begin{eqnarray}
\label{eq:2-euler}
p = (\gamma -1) (E -\frac{\rho u^2}{2}) \;,
\end{eqnarray}
for an ideal gas $\gamma = 1.4$.  Notice that, the choice of the initial condition given by the functions 
\begin{eqnarray}
\begin{array}{lcl}
\rho(x,0) &=& 1+0.2\sin(2\pi x)\;,\\
   u(x,0) &=& 1,\\
   p(x,0) &=& 2,
\end{array}
\end{eqnarray}
provides the exact solution for the system (\ref{eq:1-euler}), which corresponds to the set of functions 
\begin{eqnarray}
\begin{array}{ccl}
\rho(x,t) &=& 1+0.2\sin(2\pi (x-t))\;,\\
u(x,t)    &=& 1\;,\\
p(x,t)    &=& 2\;.
\end{array}
\end{eqnarray}
Notice that, the variables $\rho, u ,p$ correspond to the non-conservative variables, the corresponding translation to conserved variables needs to be done.  This test has a complex eigenstructure, which is a challenge for numerical methods. 
Table \ref{Table-Euler}, shows the results of the empirical convergence rate assessment for the density variable $\rho$, at $t_{out} = 1$, $\alpha = 2$ and $C_{cfl}= 0.1$, we observe that the scheme achieves the expected theoretical orders of accuracy.
\begin{table}
\begin{center}
Theoretical order : 2 \\
\begin{tabular}{cccccccc} 
\\
\hline
\hline 
Mesh  & $L_\infty$ - err & $L_\infty$- ord  & $L_1$ - err & $L_1$ - ord & $L_2$ - err & $L_2$ - ord & CPU  \\  
\hline

    16  &  $   7.03e-0 2  $  &   1.34  &  $   4.24e-0 2  $  &   1.56  &  $   4.80e-0 2  $  &   1.49  &  0.49  \\
    32  &  $   1.50e-0 2  $  &   2.23  &  $   5.53e-0 3  $  &   2.94  &  $   7.65e-0 3  $  &   2.65  &  2.50  \\
    64  &  $   3.96e-0 3  $  &   1.92  &  $   1.10e-0 3  $  &   2.33  &  $   1.57e-0 3  $  &   2.28  &  9.13  \\
   128  &  $   1.37e-0 3  $  &   1.53  &  $   3.05e-0 4  $  &   1.85  &  $   5.08e-0 4  $  &   1.63  &  39.15  \\

 \hline
 \\
 \end{tabular}  
\\
Theoretical order : 3 \\
\begin{tabular}{cccccccc}  
\\
\hline 
Mesh  & $L_\infty$ - err & $L_\infty$- ord  & $L_1$ - err & $L_1$ - ord & $L_2$ - err & $L_2$ - ord & CPU  \\  
\hline

    16  &  $   3.86e-0 2  $  &   1.90  &  $   2.52e-0 2  $  &   2.01  &  $   2.78e-0 2  $  &   1.98  &  2.77  \\
    32  &  $   5.54e-0 3  $  &   2.80  &  $   3.55e-0 3  $  &   2.83  &  $   3.94e-0 3  $  &   2.82  &  12.28  \\
    64  &  $   7.08e-0 4  $  &   2.97  &  $   4.52e-0 4  $  &   2.97  &  $   5.01e-0 4  $  &   2.97  &  46.66  \\
   128  &  $   8.89e-0 5  $  &   2.99  &  $   5.66e-0 5  $  &   3.00  &  $   6.29e-0 5  $  &   3.00  &  179.21  \\

 \hline
 \\
 \end{tabular}  
\\
Theoretical order : 4 \\
\begin{tabular}{cccccccc}  
\\
\hline 
Mesh  & $L_\infty$ - err & $L_\infty$- ord  & $L_1$ - err & $L_1$ - ord & $L_2$ - err & $L_2$ - ord & CPU  \\  
\hline

    16  &  $   9.95e-0 3  $  &   3.73  &  $   6.45e-0 3  $  &   3.85  &  $   7.12e-0 3  $  &   3.82  &  35.56  \\
    32  &  $   3.43e-0 4  $  &   4.86  &  $   2.18e-0 4  $  &   4.89  &  $   2.43e-0 4  $  &   4.87  &  154.7`  \\
    64  &  $   1.10e-0 5  $  &   4.96  &  $   7.00e-0 6  $  &   4.96  &  $   7.78e-0 6  $  &   4.96  &  354.59  \\
   128  &  $   3.58e-0 7  $  &   4.94  &  $   2.28e-0 7  $  &   4.94  &  $   2.53e-0 7  $  &   4.94  &  1397.81  \\

\hline
 \\
 \end{tabular}  
\\
Theoretical order : 5 \\
\begin{tabular}{cccccccc}  
\\
\hline 
Mesh  & $L_\infty$ - err & $L_\infty$- ord  & $L_1$ - err & $L_1$ - ord & $L_2$ - err & $L_2$ - ord & CPU  \\  
\hline

    16  &  $   1.32e-0 3  $  &   4.59  &  $   8.58e-0 4  $  &   4.69  &  $   9.47e-0 4  $  &   4.67  &   222.60  \\
    32  &  $   4.32e-0 5  $  &   4.93  &  $   2.77e-0 5  $  &   4.95  &  $   3.07e-0 5  $  &   4.95  &   839.69  \\
    64  &  $   1.36e-0 6  $  &   4.99  &  $   8.68e-0 7  $  &   4.99  &  $   9.64e-0 7  $  &   4.99  &  2500.04  \\
   128  &  $   4.21e-0 8  $  &   5.02  &  $   2.68e-0 8  $  &   5.02  &  $   2.98e-0 8  $  &   5.02  &  7757.30 \\

\hline

\end{tabular} 
\end{center}
\caption{Euler equations.Output time
  $t_{out} = 1$ with  $C_{cfl}= 0.1.$, $\alpha = 2$.}\label{Table-Euler}
\end{table}

\section{Conclusions}\label{sec:conclusions}
In this work, a new scheme of the family of the ADER methods has been presented. The scheme has been proved to achieve the expected theoretical order of accuracy up to the fifth order of accuracy and able to deal with stiff source terms. The scheme has been based on the implicit Taylor series expansion proposed in \cite{Toro:2015a} where implicit evolution for the data and its derivatives are needed. Here, the evolution of the data has been done as in \cite{Toro:2015a}, where conventional Cauchy-Kowalewsky is implemented. However, the evolution of the spatial derivatives here has been done by a second order implicit Taylor series, where the leading term is a linearization around the resolved GRP state.  The strategy has resulted in a very fast inversion matrix procedure.
 
A von Neumann analysis has been carried out, where a reduction in CFL has been observed for stability purposes. However, as shown here the improvement in the predictor stage allows an enhancement on the efficiency as the order of accuracy increases.  In order to reduce the introduction of numerical dissipation due to small CFL, the FORCE$-\alpha$ scheme has been implemented in a non-conservative fashion and thus present scheme is also a universal one. 

Despite, the values for choosing $\alpha$ for the $FORCE-\alpha$ has been reported in \cite{Toro:2020a}, there is no evidence that for high-order the same strategies are still valid. In this work, the values have been taken by the hand and further research is required.

\section*{Acknowledgments}
G.M thanks to the {\it National chilean Fund for Scientific and Technological Development}, FONDECYT, in the frame of the research project for Initiation in Research, number 11180926.

\bibliographystyle{plain}
\bibliography{ref}

\appendix

\section{The reconstruction procedure}\label{sec:reconstruction}

In this work the von Neumann Analysis of the high order scheme, is carried out. Since a critical element on this analysis regards the reconstruction procedure, we are going to review the WENO reconstruction reported in \cite{Kaser:2007a}. Furthermore, the procedure requires the use of coefficients, randomly distributed which are going to be detailed in this appendix also. To describe the process, we are going to assume that the scalar cell averages $q_i^n$ are available at each cell at time $t^n$. Then, the reconstruction polynomial $ \tilde{ w}_i(x) $ defined on $ [x_{i-\frac{1}{2}}, x_{i+\frac{1}{2} } ] $ is expressed as 
$ \tilde{w}_i(x) = \omega_L \cdot \tilde{p}_L(x) +  \omega_C \cdot \tilde{p}_C(x)  + \omega_R \cdot \tilde{p}_R(x)$,
where  $\omega_L$, $\omega_C$ and $\omega_R$ are coefficients, which depend on the data and satisfy $ \omega_L + \omega_C + \omega_R = 1$. The functions $\tilde{p}_L(x)$, $\tilde{p}_C(x)$  and $\tilde{p}_R(x)$ are polynomials defined on $[x_{i-\frac{1}{2}}, x_{1+\frac{1}{2}}]$ and constructed from sets of cells averages refereed to as stencils, as described below.
By simplicity, let us transform the interval $ [x_{i-\frac{1}{2}}, x_{i+\frac{1}{2}}]$ into $[0,1]$. This is done by the change of variable $ x = x_{i-\frac{1}{2}} + \xi \Delta x$. Therefore, the reconstruction polynomial $ w_i (\xi ) = \tilde{ w}_i(x(\xi))$ takes the form
\begin{eqnarray}
\begin{array}{c}

w_i( \xi) = \omega_L \cdot p_{L}(x) + \omega_C \cdot p_{C} (\xi) + \omega_R \cdot p_{R} (\xi)  \;,
\end{array}
\end{eqnarray}
where $ p_{S} (\xi) = \tilde{p}_S(x(\xi))$, $S\in \{ L,C,R\}$. 

The polynomial $ p_L(\xi)$ is obtained through the set of states $q^n_j$ with $j \in \Omega_L :=  \{ i-M, i- M + 1,...,i-1, i \}$, where $M+1$ is the order of accuracy. 
The polynomial $ p_R(\xi)$ is obtained through the set of states $q^n_j$ with $j \in  \Omega_R :=   \{ i, i + 1,...,i+M \}$. 

The polynomial $ p_C(\xi)$ is obtained through the set of states $q^n_j$ with $j \in  \Omega_C$. The set of indices for $M$ odd we take $\Omega_C :=   \{ i-M, i- M + 1,...i-1, i, i + 1,...,i+M \}$ and for $M $ even $\Omega_C :=   \{ i-\frac{M}{2}, i- \frac{M}{2} + 1,...i-1, i, i + 1,...,i+\frac{M}{2} \}$    Each one of these polynomials has the form
\begin{eqnarray}
\begin{array}{c}
\displaystyle
p_S (\xi) = 
\sum_{l=0}^{M} \beta^S_l \cdot \theta_l (\xi) \;,
\end{array}
\end{eqnarray}
where $\theta_l ( \xi)$ is the $lth$ Legendre polynomial on $[0,1]$, defined by
\begin{eqnarray}
\begin{array}{c}
\displaystyle

\theta_l ( \xi ) = 
(-1)^l
\sum_{k=0}^l 

\left(
\begin{array}{c}
l \\
k
\end{array}
\right)
\cdot 

\left(
\begin{array}{c}
l + k \\
k
\end{array}
\right) (- \xi )^k \;.
\end{array}
\end{eqnarray}
Therefore, the coefficients $\beta^S_l$ are found by solving the following  
\begin{eqnarray}
\begin{array}{c}
\displaystyle
q^n_j =   \sum_{l = 0}^M  \beta_l^S\int_{j-i}^{j-i+1} \theta_l(\xi)d\xi 
\;,
\end{array}
\end{eqnarray}
with $ j \in \Omega_S$. Notice that the previous expression is reduced to solve a linear system for the coefficients $ \beta^S_j$. However, for the special case of $M$ odd, the system associated to $ p_C$ is a overestimate one which is solved as a constraint least square problem, see  \cite{Kaser:2007a} for further details. 
 
Once the polynomials have been computed, the weights $ \omega_S$ are computed as 
\begin{eqnarray}
\begin{array}{c}
\displaystyle
\tilde{\omega}_S = \frac{ \lambda_S}{ ( OI_S + \varepsilon)^r} \;,  
\end{array}
\end{eqnarray}
where $r = 4$, $\varepsilon = 10^{-14}$, $ \lambda_L = \lambda_R$ and $ \lambda_C = 10^5$ and $OI_S$ is an oscillation index which is computed as 
\begin{eqnarray}
\begin{array}{c}
\displaystyle
OI_S = \sum_{k = 1}^M\int_0^1 (\frac{d^k}{d\xi^k} p_S(\xi))^2 d\xi \;.
\end{array}
\end{eqnarray}
Then $ \omega_S = \frac{\tilde{\omega}_S}{ \tilde{\omega}_L + \tilde{\omega}_C + \tilde{\omega}_R }$.

As can be seen, the coefficients $ \tilde{\omega}_S$, with $S=L,C,R$ depend on several factors. However, once a normalization is carried out to get $\omega_S$, they become a convex combination.  Therefore for the purpose of the von Neumann analysis we are going to assume that these values are randomly distributed, but with the constraints $ \omega_S \geq 0 $ and $ \omega_L + \omega_C + \omega_R= 1$.

Notice that the reconstruction polynomials $w_i(\xi)$ can be written as
\begin{eqnarray}
\begin{array}{c}
\displaystyle
w_i(\xi) 
= \sum_{k=0}^M
(\omega_L \beta_k^L + \omega_C \beta_k^C + \omega_R \beta_k^R   ) p_k(\xi) \;,

\end{array}
\end{eqnarray}
Since $\beta^S_k$ for any $S=L,C,R$ depends on the data, $q_{i+j}^n $ with $j \in \Omega_S$, the reconstruction polynomial can  be written  as
\begin{eqnarray}
\begin{array}{c}
\displaystyle
w_i( \xi) 
= \sum_{u=-M}^M
\phi_{i+u} (\xi)
\cdot q_{i+u}^n
\;.
\end{array}
\end{eqnarray}
Notice that the shape of functions $\phi_{i+u}$ is different for any order of accuracy, specifically, this depends on the
$\beta^S_j$, with $j=0,...,M$ and $S = L,C,R$. For instance in the third order case, the coefficients for the left stencil have the form
\begin{eqnarray}
\begin{array}{l}
\beta^L_0 = q_i^n\;, \\
\beta^L_1 = \frac{3}{ 4} \cdot q_{i}^n - q_{i-1}^n - \frac{1}{4} \cdot q_{i-2}^n \;, \\
\beta^L_2 = \frac{1}{12} \cdot q_{i}^n - \frac{1}{6} \cdot q_{i-1}^n + \frac{1}{12} \cdot q_{i-2}^n \;. \\
\end{array}
\end{eqnarray}
For the central stencil the coefficients have the form
\begin{eqnarray}
\begin{array}{l}
\beta^C_0 = q_i^n\;, \\
\beta^C_1 = - \frac{1}{4} \cdot q_{i-1}^n+\frac{1}{4} \cdot q_{i+1}^n \;, \\
\beta^C_2 = \frac{1}{12} \cdot q_{i-1}^n - \frac{1}{6} \cdot q_{i}^n + \frac{1}{12} \cdot q_{i+1}^n \;. \\
\end{array}
\end{eqnarray}
For the right stencil the coefficients have the form
\begin{eqnarray}
\begin{array}{l}
\beta^R_0 = q_i^n\;, \\
\beta^R_1 = - \frac{3}{4} \cdot q_{i}^n + \cdot q_{i+1}^n -  \frac{1}{4} q_{i+2}^n  \;, \\
\beta^R_2 = \frac{1}{12} \cdot q_{i}^n - \frac{1}{6} \cdot q_{i+1}^n + \frac{1}{12} \cdot q_{i+2}^n \;. \\
\end{array}
\end{eqnarray}
Collecting this information and considering that for third order we have $M=2$, the functions $\phi_{i+u}(\xi)$ correspond to
\begin{eqnarray}
\begin{array}{l}

\phi_{i-2} ( \xi) = \frac{ \omega_L}{4} \cdot ( \frac{p_2(\xi)}{3} - p_1(\xi) ) \;, \\
\\
\phi_{i-1} ( \xi) = p_1(\xi) \cdot ( - \omega_L - \frac{\omega_C}{4} ) 
           - p_2(\xi)\frac{ 1 }{6} \cdot ( \frac{\omega_C}{2} - \omega_L )\;, \\ 
\\
\phi_{i} ( \xi) = 1 + p_1(\xi) \cdot \frac{3}{4} \cdot (  \omega_L - \omega_R  ) 
                + p_2(\xi) \cdot \frac{1}{6} \cdot ( \frac{\omega_L}{2} - \omega_C + \omega_R ) \;,
\\
\\
\phi_{i+1} ( \xi) = p_1(\xi) \cdot (  \omega_R + \frac{\omega_C}{4} ) 
           + p_2(\xi)\frac{ 1 }{6} \cdot ( \frac{\omega_C}{2} - \omega_R )\;, \\ 
\\
\phi_{i+2} ( \xi) = \frac{ \omega_R}{4} \cdot ( \frac{p_2(\xi)}{3} - p_1(\xi) ) \;. \\

\end{array}
\end{eqnarray}
Since the $\omega_S$ are randomly distributed, they do not depend on the data and thus $\phi_{i+u}$ neither.

\section{A practical Cauchy-Kowalewskaya functional generator}\label{sec:ck-generator}

Any algebraic software manipulator can be implemented for obtaining the Cauchy-Kowalewskaya functions and the gradient of the non-linear algebraic system resulting from the implicit Taylor series expansion.  However, for the sake of completeness, we provided a script, that has been written to generate each one of the expressions resulting from the Cauchy-Kowalewskaya procedure required in this paper.

Maxima  \cite{Ochsner:2019a}, is a GNU computer algebra system which can be obtained at http://maxima.sourceforge.net and it is compatible with any operative system.

\begin{verbatim}

/* 
-------------------------------------------------------

            Cauchy-Kowalewskaya procedure 
                     
                         for 
                      
             a general hyperbolic system
         
-------------------------------------------------------

This script provides the Cauchy-Kowalewskaya functionals
for a partial differential equation of the form

                 Qt + A(Q) * Qx = S(Q)

The input of this script corresponds to the expression
for the advective term, "A(Q) * Qx", the source term "S(Q)"
and the order of accuracy. Therefore, the user must modify 
the following items:

- order of accuracy "N", 
- the number of unknowns, "nVar",  
- the Jacobian matrix "A",
- the source function "S".

In this particular case, we implement the procedure for
the non-conservative system.

To execute use CTR + R.


-------------------------------------------------------
*/
kill ( all);
/*  
-------------------------------------------------------
Order of accuracy: N,  N+1  degrees of freedom.
(In this example N = 3 and so the third order is set )
-------------------------------------------------------
 */
N      :  3 $
/*  
-------------------------------------------------------
Number of unknowns: nVar.  
-------------------------------------------------------
*/ 
nVar : 2 $


/*  
-------------------------------------------------------
Load package to print out into a Fortran 90 format.   
(comment this line if a .f90 format is not required)
-------------------------------------------------------
*/
load("f90");


/*  
-------------------------------------------------------
Variables (space-time evaluated) (nVar).
This part has to be modified by the user.  
In this application, two variables are taken into account
For each variable, it is needed to be included as follows:
-------------------------------------------------------
*/
u [ 1] : q1 ( x, t) $
u [ 2] : q2 ( x, t) $

/*  
-------------------------------------------------------
Jacobian  matrix A: 
This part has to be modified by the user.  
Put here in a component wise the expression A * dxQ,
which corresponds to the convective part.  
The variables are represented by "u".
-------------------------------------------------------
*/ 

A : matrix(  
           [ lambda,   u[1]   ], 
           [      1,   lambda ]  ) $

ADx[1]: A[1,1] * diff(u[1], x) + A[1,2] * diff(u[2], x) $
ADx[2]: A[2,1] * diff(u[1], x) + A[2,2] * diff(u[2], x) $

/*  
-------------------------------------------------------
Source s: 
This part has to be modified by the user.  
Put here the expression for the source function. 
The variables are represented by "u".
-------------------------------------------------------
*/
s[1]:     R * u[1] * ( u[1] - 1 ) $
s[2]: - R      *      ( u[2] - 1 ) $

/*
-------------------------------------------------------
Definition of the variable "q", which is the variable 
used to work.
This represents the time derivative of the variable "u".
"q[k,l]", k the number of variables, "l" the order of 
derivative.
Here "q[k,0] = u[k] " is the variable.
-------------------------------------------------------
*/
for k : 1 thru nVar do (
    q [ k, 0] : u [ k]
 ) $
  
/* 
-------------------------------------------------------
Computing high-order time-derivatives.
-------------------------------------------------------
*/
for k : 1 thru nVar do (
    for i : 2 thru N + 3 do (
        q [ k, i] : diff ( q [ k, i - 1], t),
        for j : i-1 thru 1 step -1 do (
            /*
            ------------------------------------------- 
            Replace time-time derivatives previously computed.
            ------------------------------------------- 
            */
            for k1 : 1 thru nVar do (   
            
                q [ k, i] : subst (  
                           [ 
                           diff ( u [ k1], t, 1, x, j) 
                           = 
                           diff ( q [ k1, 1], x, j) 
                           ],                            
                           q [ k, i] )     
                 )
            ),
            for k1 :1 thru nVar  do ( 
                q [ k, i] : subst (  
                            [ 
                            diff ( u [ k1], t, 1) 
                            = 
                            q [ k1, 1] 
                            ], 
                            q [ k, i] )
            ),
            
            q [ k, i] : ratsimp ( q [ k, i] )
        )
    ) $
 
 
 /* 
-------------------------------------------------------
Definition of the Cauchy-Kowalewskaya functionals. 
-------------------------------------------------------
*/

for k : 1 thru nVar do (
      for i : 1 thru N+1 do (
           G [ i, k] :q[ k, i - 1] 
       )
 ) $
 
 
/*
-------------------------------------------------------
 D_x_Q [ i, m ]: (i-1)-th spatial derivative of the m-th 
 component.
  
 D_t_Q  [ 1, m] : q [ m, i] :  output Fortran format for 
 the i-th time-derivative of the m-th component, in terms 
 of spatial-derivatives.
-------------------------------------------------------
*/

for d_t : 1 thru N do (
    for m : 1 thru nVar do (
        DtQ [ d_t, m] : q [ m, d_t - 1 ],
        for i : N + 1 thru 0 step -1 do ( 
            for m1 : 1 thru nVar do (
                DtQ [ d_t, m] : subst ( 
                       [ diff ( u [ m1] , x, i)   
                       = DxQ [ i+1, m1] ], DtQ [ d_t, m] ) 
                            
                   )
            )

        )
) $ 

/* 
-------------------------------------------------------
DxQ ( i, k) is the (i-1)-th space-derivative of the k-th 
component.   
DTQ(i,k) is the (i+1)th time-derivative in terms of space 
derivatives. 

To store the conventional Cauchy-Kowaleswaya functional 
use the following two command lines, this will generate 
a "txt" external file where the expressions are provided 
in "f90"  format.

(comment the line containing "f90" if .f90 format is not 
required).

-------------------------------------------------------
*/

Dt : genmatrix ( DtQ, N, nVar)                        $

with_stdout ( "DT.txt", f90( 'DTQ = float ( Dt) ) )   $

/*
-------------------------------------------------------
From the time derivatives we compose the implicit in time, 
Taylor series expansion. This generates a non-linear 
equation for the resolved state. To solve this equation 
we use a fixed point iteration procedure.
-------------------------------------------------------
*/

for m : 1 thru nVar do ( 
    FunHighOrder [m] : Q[m] - DxQ[1,m]  
                     + subst( [ DxQ[1, 1] = Q[ 1] , 
                                DxQ[1, 2] = Q[ 2] , 
                                DxQ[1, 3] = Q[ 3]  ], 
                       sum(  ( - t)^(k-1) / prod ( l, l, 1, k-1)
                       * DtQ[ k, m], k, 2, N  )  ) 
    );
    
indJ [ i, j] := diff ( FunHighOrder[ i], Q[ j] )  ;

/* 
-------------------------------------------------------
The Jacobian of the non-linear system derived from the 
implicit Taylor series expansion. If the f90 format
is not required, just comment the following two lines.
-------------------------------------------------------
 */
JacDt : genmatrix(  indJ, nVar, nVar) $
with_stdout ( "JacDT.txt", f90( 'JacDTQ = float ( JacDt) ) )  $
/* 
-------------------------------------------------------
The information is stores in an external "txt" file
containing the "f90" format.  In Linux distribution
the file is stored into the same folder where this
script is called.
-------------------------------------------------------
 */
      
\end{verbatim}

\end{document}